\documentclass[final]{siamltex}

\usepackage{amsmath,makeidx,amssymb,amscd,amsfonts}
\usepackage{graphicx,epsfig,latexsym,bbm,index}
\usepackage[active]{srcltx}
\usepackage[below]{placeins}
\def\bpn{\bigskip\par\noindent}

\def\ve{{\varepsilon}}

\def\ol{\overline}
\def\sgn{\mbox{sgn}}

\def\F{{\cal F}}
\def\Fe{{\cal F}_{\ve,\alpha}}
\def\Ga{{\cal G}_\alpha}
\def\Ger{{\cal G}_{\ve,\alpha}}

\def\cp{{\cal{P}}}
\def\.{{\;}}
\def\bv{{\tt{BV}}}
\def\lzo{{L^2(\Omega)}}

\def\pie{{\phi_{\ve,\alpha}}}
\def\piek{{\phi_{\ve_k,\alpha_k}}}

\def\req#1{{\rm(\ref{eq:#1})}}

\newcommand\R{\mathbb{R}}

\newcommand\N{\mathbb{N}}

\newtheorem{example}[theorem]{Example}

\newtheorem{remark}[theorem]{Remark}
\newtheorem{conjecture}[theorem]{Conjecture}

\title{\bf Analysis of Regularization Methods for the Solution of
           Ill--Posed Problems involving Discontinuous Operators}

\author{F. Fr\"uhauf, O. Scherzer\thanks{Department of Computer Science,
        Universit\"{a}t Innsbruck,
        Technikerstra{\ss}e 25, A-6020 Innsbruck, Austria,
        {\tt \{Florian.Fr\"uhauf,Otmar.Scherzer\}@uibk.ac.at}}
        \and A.~Leit\~ao\thanks{
        Radon Institute for Comput. Appl. Math., Austrian
        Academy of Sciences, c/o Kepler Universit\"{a}t Linz,
        A-4040 Linz, Austria, {\tt Antonio.Leitao@oeaw.ac.at}  }}

\date{today}

\begin{document}

\maketitle

\begin{abstract}
We consider a regularization concept for the solution of ill--posed
operator equations, where the operator is composed of a continuous
and a discontinuous operator.
A particular application is level set regularization, where
we develop a novel concept of minimizers.
The proposed level set regularization is capable of handling changing topologies.
A functional analytic framework explaining the splitting of topologies
is given. The asymptotic limit of the level set regularization method is an
evolution process, which is implemented numerically and the quality
of the proposed algorithm is demonstrated by solving an inverse source problem.
\end{abstract}



\pagestyle{myheadings}
\thispagestyle{plain}
\markboth{F. Fr\"uhauf, O. Scherzer, A. Leit\~ao}
{Regularization methods for ill--posed problems involving discontinuous
operators}

\section{Introduction} \label{sec:introd}

The goal of this paper is to analyze \emph{regularization} models for the
\emph{stable} solution of \emph{ill-posed} operator equations
\begin{equation} \label{eq:op-u}
F(P(\phi)) = y\;.
\end{equation}
Here $F$ is a continuous operator between Banach spaces $X$ and $Y$ and $P$
is a probably \emph{discontinuous operator} from an admissible
class ${\cal{P}}$ into $X$.
Classical results on convergence and stability of variational
regularization principles for solving \emph{non-linear} ill--posed
problems (see e.g. \cite{Mor84,Mor93,EngHanNeu96}) in a Hilbert spaces
setting such as
\begin{remunerate}
\item existence of a regularized solution,
\item stability of the regularized approximations,
\item approximation properties of the regularized solutions
\end{remunerate}
are applicable if the operator $P$ is
\begin{remunerate}
\item bounded and linear or
\item nonlinear, continuous, and weakly closed.
\end{remunerate}
In this paper we particularly emphasize on operator equations \req{op-u} where the
operator $P$ is discontinuous. Of particular interest for this paper is
\begin{equation}
\label{eq:pr}
P(t) := \left\{ \begin{array}{rcl}
0 & \mbox{ for } & t < 0\,,\\
1 & \mbox{ for } & t \geq 0\;.
\end{array} \right.
\end{equation}
With $P$ there is associated the \emph{admissible class}
\begin{equation}
\label{eq:cp}
\cp := \left\{ u : u =
               \chi_D \mbox{ where }
               D \subseteq \Omega \mbox{ is measurable and }
               {\cal{H}}^{n-1}(\partial D) < \infty
       \right\}\;.
\end{equation}
Here
\begin{remunerate}
\item ${\cal{H}}^{n-1}(\partial D)$ denotes the
      $n-1$-dimensional Hausdorff-measure of the boundary of
      $D$;
\item $\chi_D$ denotes the characteristic function of the set
      $D$.
\end{remunerate}
We call a regularization approach involving this projection
\emph{level set regularization} since we recover the boundary of
an object $\partial D$, which is a level set (for instance with
value $0$) of a function $\phi$. The idea of considering
characteristic functions as level sets of higher dimensional data
has been used before in the context of \emph{multiphase flow} (see
e.g. \cite{MerBenOsh94,ZhaChaMerOsh96,CheMerOshSme97}) and
segmentation (see e.g. \cite{ChaSheVes03}). Level set method have
been used successively in many applications since the pioneering
work of Osher \& Sethian \cite{OshSet88}. For solving inverse
problems applications with level sets we refer to Santosa
\cite{San96} and Burger \cite{Bur01}. \bpn In this work we base
our considerations on ideas from nonlinear convex semigroup theory
(cf. Brezis \cite{Bre73}) which allows to characterize the
solution of an evolution process by implicit time steps of
regularization models. Since our regularization models appear to
be nonconvex, the theoretical results of nonlinear semigroup
theory are not available. Simulating this approach, we show in
this work that iterated regularization is well-posed, and (aside
form the lack of theoretical results) we can interpret the
iterated regularized solutions as time instance of an evolution
process.
\bpn
Various other models fit in the general framework of this paper
but are not particularly emphasized: For instance for $a \in \R$ let us
consider the following projection operator
\[
P^a(t) := \left\{ \begin{array}{rcl}
-a & \mbox{ for } & t < -a\,,\\
t & \mbox{ for } & -a \leq t \leq a\,,\\
a & \mbox{ for } & t > a\,,
\end{array} \right.
\]
with the admissible class
\begin{equation}
\label{eq:cp1} \cp_a := \left\{ u : u = P^a(\phi) \mbox{ with }
\phi \in H^1(\Omega) \right\}\;.
\end{equation}
The operator $P^a$ ensures that the recovered functions are absolutely
bounded by $a$.
\[
P^+(t) := \exp (t)
\]
with the \emph{admissible class}
\begin{equation}
\label{eq:cp2} \cp_+ := \left\{ u : 0 < u = P^+(\phi) \right\}
\end{equation}
can be used to guarantee \emph{non-negativity}.
Depending of the operator $P$ we actually solve a constraint optimization problem.
With $P_+,P_a,P$ we guarantee that the solution is in the according admissible class.
\bpn
The outline of this paper is as follows: In Section \ref{sec:lsr}
we introduce the concept of \emph{level set regularization}, based
on considerations in \cite{San96,Bur01,LeiSch03}.
The level set regularization functionals derived in \cite{LeiSch03} are
modified such that a convergence analysis becomes tractable (cf. Section \ref{sec:tow}).
That is we show that each implicit time step is well-defined. This a prerequisite step
in showing that the according evolution process is well-defined.
To this end, we introduce a novel concept of a minimizer of regularization functionals
involving discontinuous operators (cf. Section \ref{sec:novel}). A
convergence analysis is presented in Section \ref{sec:2.3}. The
problem of numerical minimization is discussed in Section \ref{sec:2.4} and
finally numerical examples are presented in Section \ref{sec:2.5}.

\section{Analysis of Level Set Regularization}
\label{sec:lsr}
In the following we pose the general assumptions which we assume
to hold all along this paper:
\begin{remunerate}
\item $\Omega \subseteq \R^n$ is bounded with $\partial \Omega$
      piecewise $C^1$ (see e.g. \cite{Ada75}).
\item The operator
      $
      F: L^1(\Omega) \to Y
      $
      is continuous and Fr\'echet-differentiable. $Y$ is a Banach
      space.
\item $\ve,\alpha,\beta$ denote positive parameters.
\item We use
      the following notation:
      \begin{romannum}
      \item $\to$ denotes strong convergence,
      \item $\stackrel{(*)}{\rightharpoonup}$ denotes weak$({}^*)$ convergence,
      \item $L^p(\Omega)$ denotes the space of measurable, $p$-times
            integrable functions,
      \item $W^{1,p}(\Omega)$ denotes the Sobolev space of one time
            weakly differentiable functions where the function and its
            derivative are in $L^p$; in particular we set $H^1=W^{1,2}$.
      \item $\bv(\Omega)$ denotes the space of functions of
            \emph{bounded variation}.
      \end{romannum}
\item \label{ass:5}
      We assume that (\ref{eq:op-u}) has a solution, i.e. there exists
      a $z\in{\cal P}$ satisfying $F(z)=y$ and a function $\phi \in
      H^1(\Omega)$ satisfying $|\nabla \phi| \neq 0$ in a neighborhood of
      $\{\phi = 0\}$ and $P(\phi)=z$.
      If $z = \chi_A$ and $\emptyset \neq A$, then we let
      $$ \phi = - d_{\ol{A}} + d_{\ol{CA}}$$
      where $d_{\ol{A}}$ and $d_{\ol{CA}}$ denote the distance functions
      from $\ol{A}$, and ${\ol{CA}}$, respectively.
      Since $d_{\ol{A}}$ and $d_{\ol{CA}}$ are uniformly Lipschitz continuous
      (see e.g. \cite{DelZol94}), they are in $L^\infty (\Omega)$.
      Moreover, $|\nabla d_{\ol{A}}| \leq 1$ and $|\nabla d_{\ol{CA}}| \leq 1$
      (see again e.g. \cite{DelZol94}). In particular this shows that
      $d_{\ol{A}},d_{\ol{CA}} \in W^{1,\infty}(\Omega) \subseteq H^1(\Omega)$.
      Thus $z \in{\cal P}$ if $z = \chi_A$ and $A$ satisfies that
      the closure of the interior of $A$ is the closure of $A$.
 \end{remunerate}
We consider the unconstrained inverse problem of solving \req{op-u} with
\[
\begin{array}{rl}
P : H^1(\Omega) &\to \cp\;.\\
     \phi &\mapsto
           \frac{1}{2} + \frac{1}{2} \sgn (\phi) =:
           \frac{1}{2} + \frac{1}{2} \left\{ \begin{array}{rl}
                                 1 \mbox{ for } \phi \geq 0\\
                                 -1 \mbox{ for } \phi < 0
                                 \end{array} \right.
\end{array}
\]
The standard form of Tikhonov regularization for solving \req{op-u}
consists in minimizing the functional
\begin{equation}
\label{eq:tik} {\cal{F}}_\alpha (\phi) := \|F(P(\phi)) -
y^\delta\|_{Y}^2 + \alpha \|\phi - \phi_0\|_{H^1(\Omega)}^2
\end{equation}
over $H^1(\Omega)$.
Actually, we understand the minimizer $\phi_\alpha$ of this functional as
\[
\phi_\alpha = \lim_{\ve \to 0+} \pie\,,
\]
where the limit is understood in an appropriate sense (weak,
weak${}^*$ convergence) and $\pie$ minimizes the functional over
$H^1(\Omega)$
\begin{equation} \label{eq:ve}
\Fe (\phi):= \|F(P_\ve(\phi)) - y^\delta\|^2_{Y} + \alpha \|\phi -
\phi_0\|_{H^1(\Omega)}^2\,,
\end{equation}
where we use
\[
P_\ve(\phi) := \left\{
\begin{array}{rcl}
0 & \mbox{ for } & \phi < -\ve\,,\\
1 + \frac{\phi}{\ve} & \mbox{ for } & \phi \in \left[-\ve,0\right]\,,\\
1 &\mbox{ for } & \phi > 0\,,
\end{array} \right.
\]
for approximating $P$ as $\ve \to 0^+$. In this case we define
\[
P'(t) := \lim_{\ve \to 0+} P'_\ve(t) = \delta(t)\;.
\]
Here and in the following $\delta(t)$ denotes the one-dimensional
$\delta$-distribution.

Taking into account that
\[
\|P_\ve(\phi_k) - P_\ve(\phi)\|_{L^{1}(\Omega)} \leq \frac{1}{\ve}
\mbox{meas} (\Omega)
\|\phi_k - \phi\|_{\lzo}\,,
\]
the proof of existence of a minimizer of the functional $\F_{\ve,\alpha}$
is similar to the proof of existence of regularized solutions of Tikhonov
functionals for approximately minimizing nonlinear ill--posed problems
in \cite{EngKunNeu89,SeiVog89} (see also \cite{EngHanNeu96}).
\begin{theorem}
\label{th:existenz-F} For any $\phi_0 \in H^1(\Omega)$ the
functional $\Fe$ (cf. \req{ve}) attains a minimizer $\pie$ in
$H^1(\Omega)$.
\end{theorem}

\subsection{Towards an Analysis of Level Set Regularization Techniques}
\label{sec:tow} In the following we outline the difficulties in
performing a rigorous analysis for the functional
${\cal{F}}_\alpha$, defined in \req{tik}.
\begin{remunerate}
\item $\pie$ satisfies,
      \[
      \|P(\pie)\|_{L^\infty} \leq 1 \mbox{ and }
      \|\pie-\phi_0\|_{H^1(\Omega)} < \infty\;.
      \]
      Since $L^\infty(\Omega)$ is the dual of $L^1(\Omega)$, i.e.,
      $L^1{}^*(\Omega)=L^\infty(\Omega)$,
      we find that there exists a subsequence $\{\piek\}_{k \in \N}$
      such that
      \[
      \piek \rightharpoonup \phi \mbox{ in } H^1(\Omega) \mbox{ and }
      P(\piek) \stackrel{*}{\rightharpoonup} z \mbox{ in } L^\infty(\Omega)\;.
      \]
      There is no analytical evidence for $z\in\cp$, i.e. it
      may not be in the range of the operator $P$.
\item To overcome this difficulty let us assume that the sequence
      $\{\piek\}_{k \in \N}$ satisfies that the Hausdorff measure of the
      boundary of the set
      \[
      \{x : \piek (x) \geq 0\}
      \]
      is uniformly bounded.
      Then the bounded variation semi-norm of $P(\piek)$ is uni\-formly
      bounded, and consequently $P(\piek)$ has a convergent subsequence in
      $L^1(\Omega)$ showing that $z$ is admissible.
\end{remunerate}
This suggests to incorporate in the functional \req{tik} as an
additional regularization term the bounded variation semi-norm of
$P (\phi)$, penalizing the length of the zero level set of $\phi$.
Actually in design problems the necessity of incorporating such a
term is well documented in
\cite{KohStr86-I,KohStr86-II,KohStr86-III}. This leads to the
following modified regularization method of minimizing
\begin{equation} \label{eq:m}
\Ga (\phi):= \|F(P(\phi)) - y^\delta\|^2_{Y} +
             2\beta \alpha |P(\phi)|_\bv +
             \alpha \|\phi - \phi_0\|^2_{H^1(\Omega)} \;.
\end{equation}
In order to guarantee existence of a minimizer of $\Ga$ we introduce a
novel concept of a minimizer:

\subsection{Minimizing Concept}
\label{sec:novel}
\bpn
\begin{definition}
\label{def:eq:m}
\begin{remunerate}
\item A {\bf pair} of functions
      \[
      (z,\phi) \in L^\infty(\Omega) \times H^1(\Omega)
      \]
      is called \emph{admissible}
      \begin{romannum}
      \item \label{ita}
            if there exists a sequence $\{\phi_k\}_{k \in \N}$ in
            $H^1(\Omega)$ such that
            $\phi_k \to \phi$ with respect to the
            $L^2(\Omega)$-norm and
      \item \label{itb}
            if there exists a sequence $\{\ve_k\}_{k \in \N}$ of positive
            numbers converging to zero such that
            \[
            P_{\ve_k} (\phi_k) \to z \mbox{ in } L^1(\Omega)\;.
        \]
      \end{romannum}
\item A minimizer of $\Ga$ is considered any admissible pair of functions
      $(z,\phi)$ minimizing
      \begin{equation}
      \label{eq:gzphi}
      {\cal{G}}_\alpha (z,\phi) =
      \|F(z) - y^\delta\|_Y^2 + \alpha \rho(z,\phi)
      \end{equation}
      over all admissible pairs. Here
      \begin{equation}
      \label{eq:phi}
      \rho(z,\phi) := \inf \liminf_{k \to \infty}
        \left\{ 2 \beta |P_{\ve_k} (\phi_k)|_\bv +
                \|\phi_k-\phi_0\|_{H^1(\Omega)}^2 \right\}\,,
      \end{equation}
      where the infimum is taken with respect to all sequences
      $\{\ve_k\}_{k \in \N}$ satisfying Item 1(ii) and
      $\{\phi_k\}_{k \in \N}$ satisfying Item 1(i).
\end{remunerate}
A generalized minimizer of ${\cal{G}}_\alpha(\phi)$ is a minimizer of
${\cal{G}}_\alpha (z,\phi)$ on the set of admissible pairs.
\end{definition}

\begin{figure}[h]
      \includegraphics[width=0.45\textwidth]{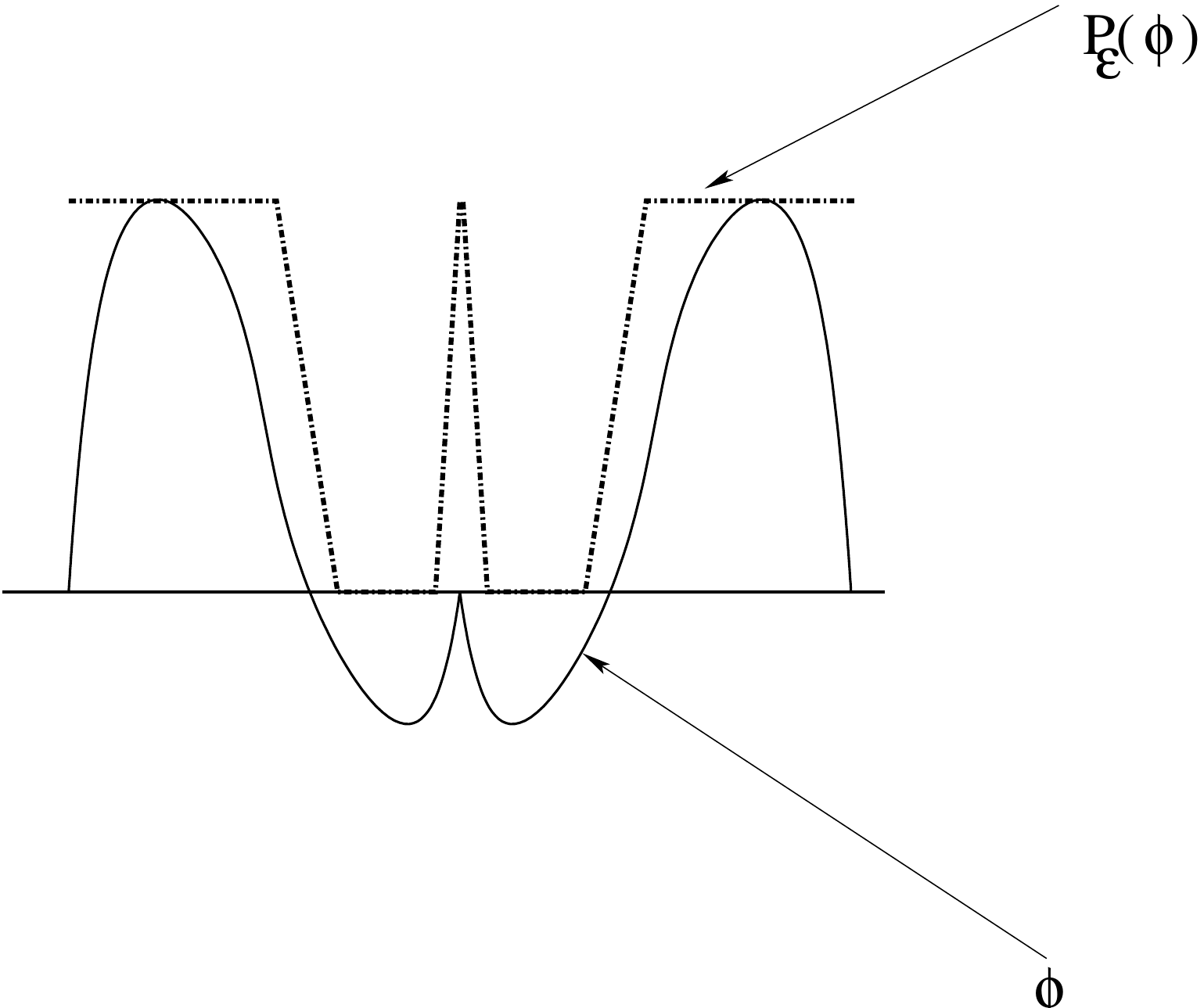} \quad
      \includegraphics[width=0.35\textwidth]{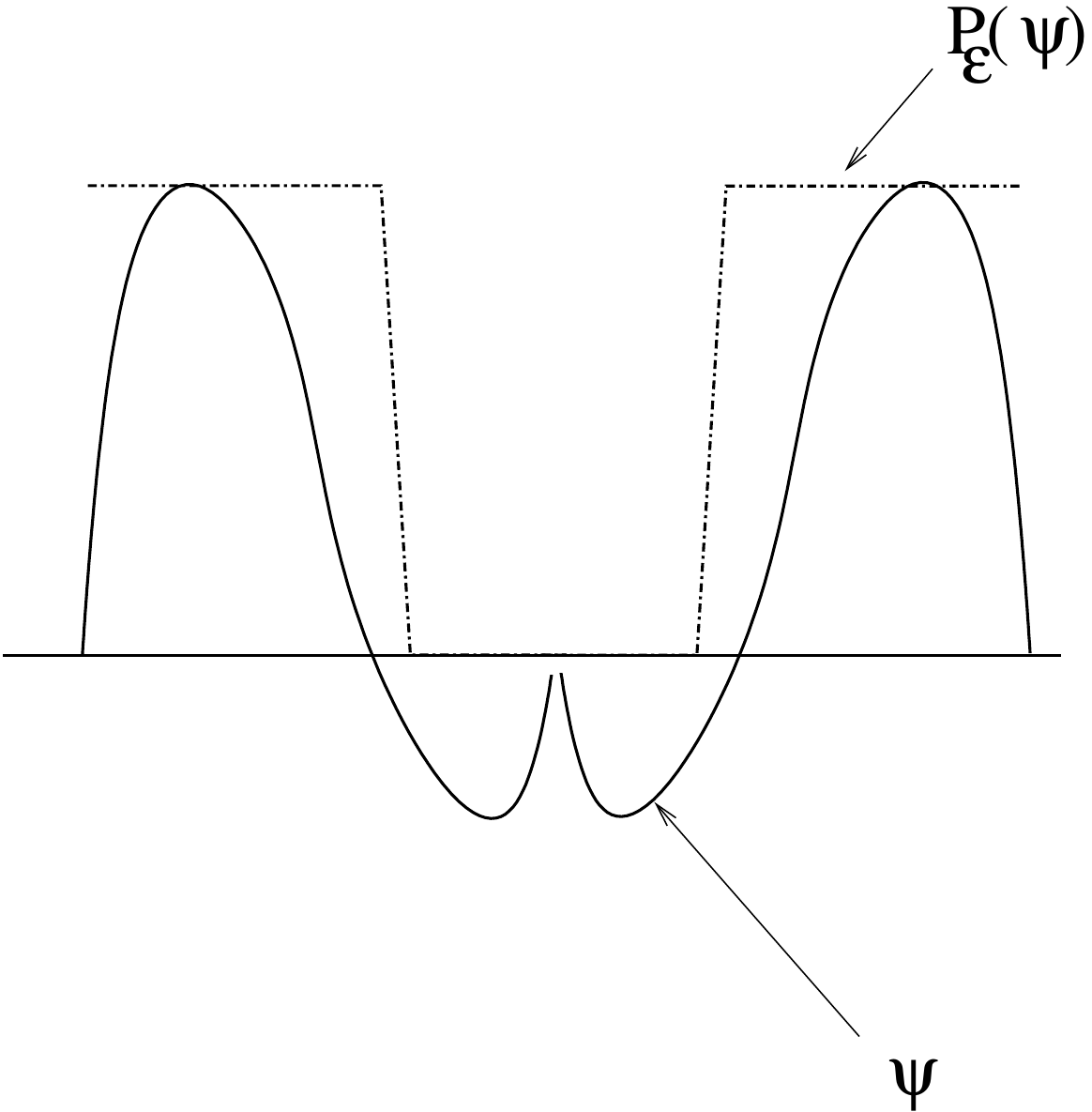}
      \caption{\label{fig:ms-1d} $n=1$:
               The functions $\phi$ and $P_\ve(\phi)$ (left):
               $|P_\ve(\phi)|_{\bv} = 4$.
               A slight perturbation: $\psi$ and
               $P_\ve(\psi)$ (right): $|P_\ve(\psi)|_{\bv}=2$.}
      \end{figure}

\begin{figure}[h]
\includegraphics[width=0.45\textwidth]{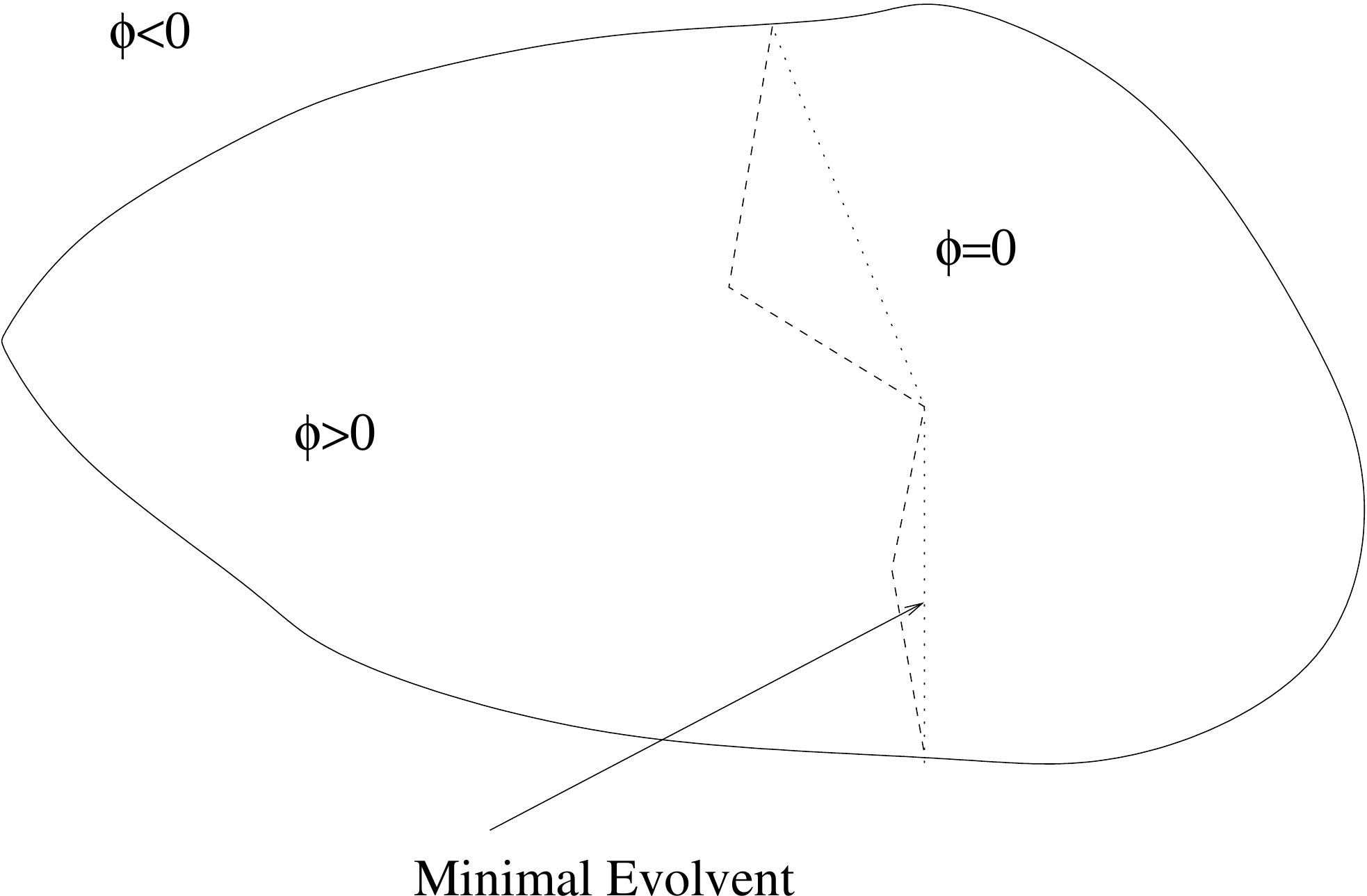} \qquad
\includegraphics[width=0.35\textwidth]{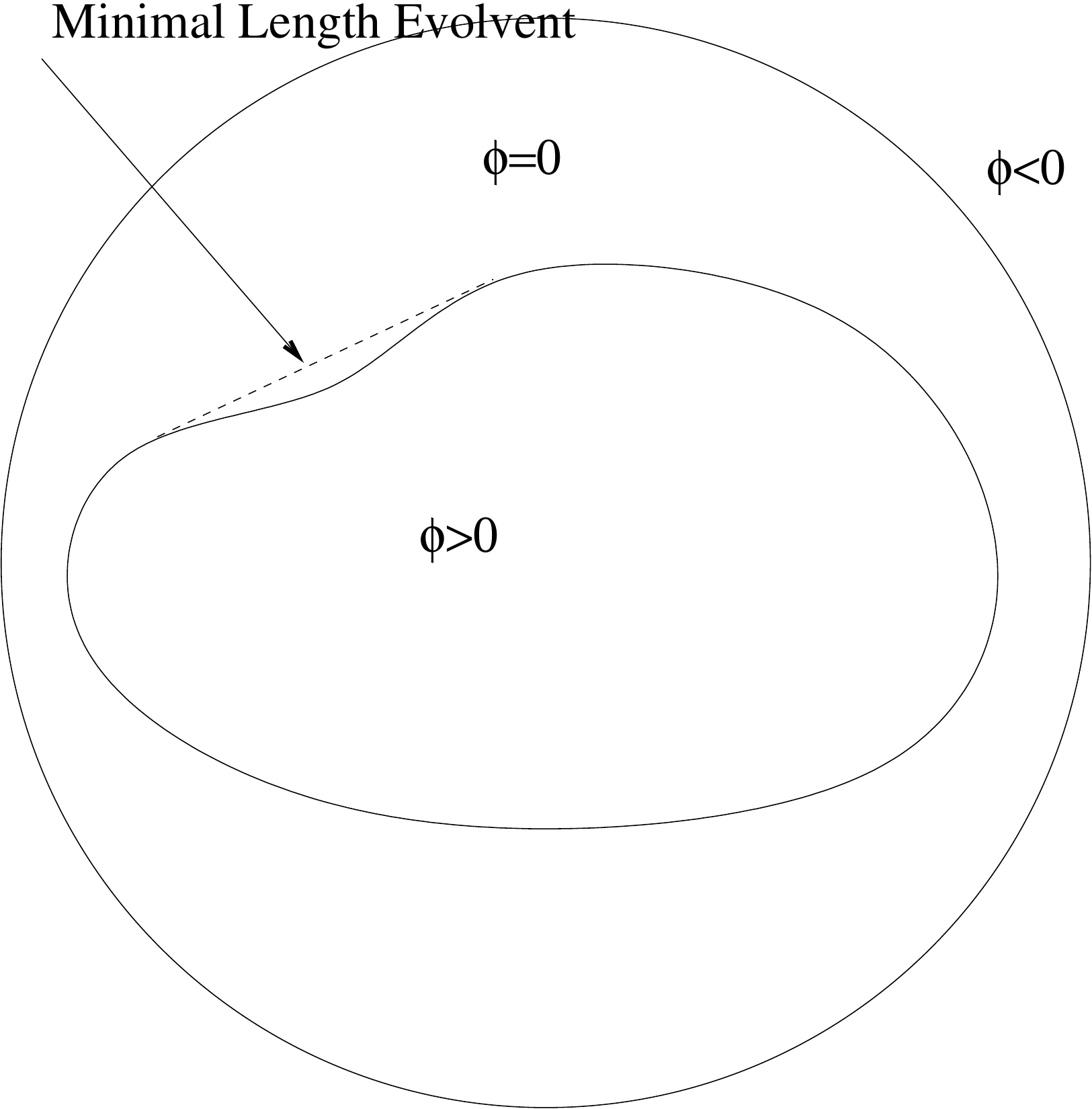}
\caption{\label{fig:ms} The minimal evolvent in $C\Phi$}
\end{figure}

The following lemma is to show that the functional $\rho$ is coercive on
the set of admissible pairs.
\begin{lemma}
\label{re:properties} For each $(z,\phi)$ admissible
\[
2\beta|z|_\bv + \|\phi - \phi_0\|^2_{H^1(\Omega)} \leq \rho(z,\phi)\;.
\]
\end{lemma}
\begin{proof}
Let $(z,\phi)$ be an admissible pair, then there exists sequences
$\{\ve_k\}_{k \in \N}$ and $\{\phi_k\}_{k \in \N}$ satisfying Items
1(i) and 1(ii) and
\[
\rho(z,\phi) = \lim_{k \to \infty}  2 \beta |P_{\ve_k} (\phi_k)|_\bv +
               \|\phi_k-\phi_0\|_{H^1(\Omega)}^2\;.
\]
By the weak lower semi-continuity of the $\bv$ and $H^1$-norms it follows
that
\[
\begin{aligned}
\| \phi - \phi_0 \| _{H^1(\Omega)}^2 & \leq
\liminf_{k \in \N} \|\phi_k-\phi_0\|_{H^1(\Omega)}^2\\
|z|_{\bv} & \leq \liminf_{k \in \N} |P_{\ve_k} (\phi_k)|_\bv\,,
\end{aligned}
\]
which proves the assertion.\qquad
\end{proof}
\bpn The definition of $\rho(z,\phi)$ is impractical, since it is
defined via a relaxation procedure. The following arguments allow
an explicit characterization of this functional. From several
experiments which we outline below, we conjecture the following
characterization of the functional $\rho (z,\phi)$.
\begin{conjecture}
      We denote by
      \[
      \Phi_+ = \{ x\in\Omega: \phi(x) > 0\} \mbox{ and }
      \Phi_- = \{ x\in\Omega: \phi(x) < 0\}
      \]
      and
      \[
      C \Phi = \Omega \backslash (\Phi_+ \cup \Phi_-)\;.
      \]
      \begin{romannum}
      \item If $\partial \Phi_+ \cap \Omega = \partial \Phi_- \cap \Omega$,
            then
            \[
            \begin{aligned}
            \rho(z,\phi) &=
                2\beta{\cal{H}}^{n-1}(\partial \Phi_- \cap \Omega) +
                \|\phi-\phi_0\|_{H^1(\Omega)}^2\\
                         &=
                2\beta {\cal{H}}^{n-1}(\partial \Phi_+ \cap \Omega) +
                \|\phi-\phi_0\|_{H^1(\Omega)}^2\;.
        \end{aligned}
        \]
      \item If the $n$-dimensional Lebesgue measure $\lambda^n(C \Phi)>0$, then $z$ is not unique
            identified, in particular $z$ can attain all values in
            $[0,1]$ in $C \Phi$. We conjecture, that
            \[
            \inf_{z\mbox{ admissible}} \rho(z,\phi)= 2\beta {\cal
            H}^{n-1}(S) + \|\phi-\phi_0\|^2_{H^1(\Omega)}.
            \]
            The problem consists in finding the surface $S$ of minimal
            $n-1$-dimensional Hausdorff measure, which is
            contained in $C\Phi$ and divides $\Omega$ in two
            sets. One set completely contains $\Phi_+$ and the other set
            contains $\Phi_-$, (cf. Figures \ref{fig:ms-1d} and
            \ref{fig:ms}).
      \end{romannum}
\end{conjecture}
Intuitively the conjecture is quite obvious. Assuming the conjecture
to be true we are further led to conjecture that the functional $\rho$ is
independent of the choice of the approximation $P_\ve$.
Thus any other approximation of $P$ with Lipschitz-continuous functions
$P_\ve$ approximating the $\delta$-distribution is suitable as well.
\begin{remark}
\label{re:25} For $\phi \in H^1(\Omega)$ where $\{\phi=0\}$ is a set
of positive Lebesgue measure (cf. Figure \ref{fig:paradox}) it is possible
to find sequences $\{\phi_k\}_{k\in\N}$ and $\{\tilde{\phi}_k\}_{k\in\N}$, which
converge strongly to $\phi$ in $\lzo$, respectively. But the limits of the
projections are different, i.e.,
$z=\lim_{k\rightarrow\infty} P_{\ve_k}(\phi_k) \neq
\tilde{z}=\lim_{k\rightarrow\infty} P_{\ve_k}(\tilde{\phi}_k)$,
cf. Figure \ref{fig:paradox}. In such a situation we have
$\rho(z,\phi) \neq \rho(\tilde{z},\phi)$.
\begin{figure}[h]
\begin{center}
\includegraphics[width=0.50\textwidth]{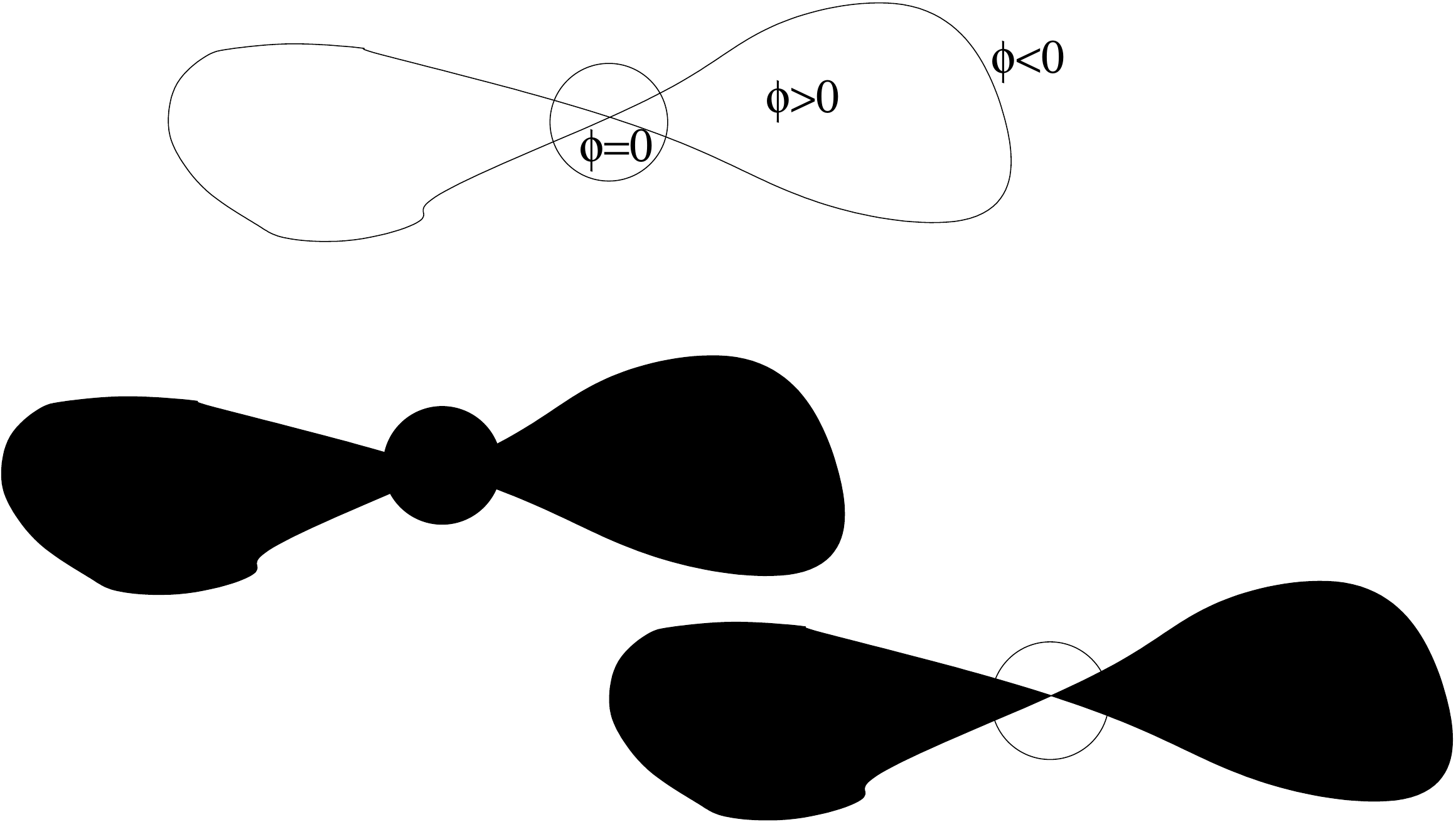} \quad
\end{center}
\caption{\label{fig:paradox} \sl{Top}: The level set function has
critical values (i.e. $|\nabla \phi|=0$ in a circle). \sl{Bottom}:
Two possible functions $z$ and $\tilde{z}$. The black value
corresponds to a value of $z = 1$.}
\end{figure}
\end{remark}

In the following we summarize some properties of the functional
$\rho$.
\begin{lemma}
\label{le:rho}
The functional $\rho$ satisfies
      \[
      \rho(z,\phi) \leq \liminf_{n \in \N} \rho (z_n,\phi_n)
      \]
      if $z_n \to z$ in $L^1(\Omega)$ and if $\phi_n \rightharpoonup \phi$ in
      $H^1(\Omega)$ and $(z_n,\phi_n)$ is admissible.
\end{lemma}
\begin{proof}
>From the definition of the functional $\rho$ and Lemma \ref{re:properties}
it follows that the functional $\rho$ is a $\Gamma^-$-limit (see e.g. \cite{Amb99})
and thus it is weak lower semi-continuous.\qquad
\end{proof}
\begin{remark}
Suppose for the moment that $P$ is a continuous operator, in which case we can
set $P_{\ve} := P$. Then the admissible class is just the set of pairs
$(z,\phi)$ satisfying $P(\phi)=z\;.$
This is just another formulation of constraint optimization.
In our context $P$ is discontinuous and therefore we consider the more
general concept of admissible pairs.
\end{remark}
\begin{example}
Let $\phi \in H^1(\Omega)$ satisfying $|\nabla \phi| > 0$ in a neighborhood
of $\{\phi=0\}$.
\begin{remunerate}
\item Let $\phi_k = \phi \in H^1(\Omega)$ and let $z = P(\phi_k)$.
      Since for any sequence $\ve_k \to 0$
      \[
      P_{\ve_k}(\phi_k) \to z\mbox{ in }L^1(\Omega)\,,
      \]
      it follows that $(z,\phi)$ is admissible.
\item Let $\phi \in H^1(\Omega)$ and denote by $\phi_k = \frac{1}{k} \phi$.
      Then there is a sequence $\ve_k \to 0$ with
      \[
      P_{\ve_k}(\phi_k) \to z\mbox{ in }L^1(\Omega)\;.
      \]
      Consequently, $(z,0)$ is admissible.
\end{remunerate}
The consequence of the second item is striking. Suppose that
$\phi_0=0$ and that there exists a minimizer $\phi_\alpha \neq 0$
of \req{m}. Then for any $k \in \N$
\[
\Ga(\phi_\alpha/k) < \Ga(\phi)\,,
\]
showing that a minimizer of $\Ga$ is not attained in a common
setting. However, the pair $(z=P(\phi_\alpha),0)$ is admissible
and can be considered as the generalized solution.

Note that in this example we consider only functions $\phi \in H^1(\Omega)$
without critical points along the zero level set.
\end{example}

\subsection{Well-Posedness and Convergence Analysis}
\label{sec:2.3}
\bpn
\begin{theorem}[{\rm Well-Posedness}]
\label{th:admissible}  Both the functional $\Ga$ and the
functional
\[
\tilde{\Ga}(z,\phi):= \|F(z)-y^\delta\|_Y^2 + 2 \beta \alpha
|z|_\bv + \alpha \|\phi-\phi_0\|_{H^1(\Omega)}^2
\]
attain minimizers on the set of admissible pairs.
\end{theorem}
\begin{proof}
\begin{remunerate}
\item Since $(0,0)$ is admissible, the set of admissible pairs is
      not empty.
\item Suppose that $\{(z_k,\phi_k)\}_{k \in \N}$ is a
      sequence of admissible pairs such that
      \[
      {\cal{G}}_\alpha(z_k,\phi_k) \to \inf {\cal{G}}_\alpha
      \leq {\cal{G}}_\alpha(0,0) < \infty\;.
      \]
      From Lemma \ref{re:properties} it follows that
      $\{(z_k,\phi_k)\}_{k \in \N}$
      is uniformly bounded in $\bv \times H^1(\Omega)$.
      By the Sobolev embedding theorem there exists a subsequence,
      denoted again by $\{\phi_k\}_{k \in \N}$ such that
      \[
      \begin{aligned}
      \phi_k \rightharpoonup \phi \mbox{ in } H^1(\Omega) \mbox{ and }
      \phi_k \to \phi \mbox{ in }L^2(\Omega)\,,\mbox{ and }\\
      z_k \to z \mbox{ in } L^1(\Omega)\,,\;
      2 \beta |z|_\bv \leq \rho(z,\phi) \leq
      \liminf_{k\rightarrow\infty} \rho(z_k,\phi_k)\;.
      \end{aligned}
      \]
      Since $\rho$ is weakly lower semi-continuous
      (cf. Lemma \ref{le:rho}) it follows that
      \begin{equation}
      \label{eq:no}
      \begin{aligned}
      \inf \Ga
      &= \lim_{k \to \infty} {\cal{G}}_\alpha(z_k,\phi_k)\\
      &= \lim_{k \to \infty}
         \left\{\|F(z_k)-y^\delta\|^2_Y
                + \alpha \rho(z_k,\phi_k) \right\}\\
      &\geq \|F(z)-y^\delta\|^2_Y +
            \alpha \rho(z,\phi)\\
      &=\Ga(z,\phi)\;.
      \end{aligned}
      \end{equation}
\item It remains to prove that $(z,\phi)$ is admissible.
      For $k$ fixed; since $(z_k,\phi_k)$ is admissible there
      exists
      a sequence $\{\ve_{k,l}\}_{l \in \N}$ of positive numbers and
      a sequence $\{\phi_{k,l}\}_{l \in \N}$ in $H^1(\Omega)$ such that
      \[
      \phi_{k,l} \to_{l \to \infty} \phi_k \mbox{ in }L^2(\Omega)\,,\quad
      P_{\ve_{k,l}} (\phi_{k,l}) \to_{l \to \infty} z_k \mbox{ in }
      L^1(\Omega)\;.
      \]
      Thus there exists an index $l(k) \in \N$ such that
      \begin{romannum}
      \item $\ve_{k,l(k)} < \frac{1}{2}\ve_{k-1,l(k-1)}$;
      \item $\|\phi_{k,l(k)} - \phi_k\|_{L^2(\Omega)} \leq \frac{1}{k}$;
      \item $\|P_{\ve_{k,l(k)}}(\phi_{k,l(k)}) - z_k\|_{L^1(\Omega)}
             \leq \frac{1}{k}$.
      \end{romannum}
      Define
      \[
      \psi_k := \phi_{k,l(k)} \mbox{ and } \eta_k := \ve_{k,l(k)}\;.
      \]
      Then, since
      \[
      \psi_k \to \phi \mbox{ in } L^2(\Omega) \mbox{ and }
      P_{\eta_k} (\psi_k) \to z \mbox{ in } L^1(\Omega)\,,
      \]
      we see that $(z,\phi)$ is admissible.
\end{remunerate}
The proof of existence of a minimizer of $\tilde{\Ga}$ is analogous as
for $\Ga$, and thus omitted.\qquad
\end{proof}
\bpn
We have shown that for any positive parameters $\alpha,\beta$ the
functionals $\Ga$ and $\tilde{\Ga}$ both attain a minimizer.
\bpn
In the sequel we denote by $(z_\alpha,\phi_\alpha)$ a minimizer of $\Ga$.

In the following we summarize some convergence result for the regularized
minimizers, which are based on the existence of a \emph{minimum norm solution}:
\begin{theorem}[{\rm Existence of a minimum norm solution}]
\label{th:min norm} Under the general assumptions of this paper
there exists a minimum norm solution $(z^\dag,\phi^\dag)$, that is
an admissible pair of functions that satisfies
\begin{remunerate}
\item $F(z^\dag )=y$,
\item $\rho(z^\dag , \phi^\dag )= \mbox{ms}:=\inf\left\{\rho(z,\phi):
      (z,\phi)\mbox{ admissible and } F(z)=y\right\}\;.$
\end{remunerate}
\end{theorem}
\begin{proof}
\begin{remunerate}
\item According to assumption \ref{ass:5} there exists a function
      $\tilde{z}\in {\cal P}$ and a function
      $\tilde{\phi}\in H^1(\Omega)$ such that $P(\tilde{\phi})=\tilde{z}$ and
      $F(\tilde{z})=y$.
      Then the pair $(\tilde{z},\tilde{\phi})$ is admissible for
      the sequence $\tilde{\phi}_k =\tilde{\phi}$, because
      $P_{\varepsilon_{k}}(\tilde{\phi}_k)\rightarrow \tilde{z}$ converges in
      $L^1(\Omega)$ for every sequence $\varepsilon_k\rightarrow 0$ due to the
      fact that $P_{\ve_k}$ is a convolution of $P$ with a $\delta$-
      distribution,
      i.e. $P_{\ve_k}=P\ast \delta_k$.
      Thus the set of admissible pairs with $F(z)=y$ is not empty.
\item Suppose that $\{(z_k,\phi_k)\}_{k\in\N}$ is a sequence of
      admissible pairs with $F(z_k)=y$  such that
      \[
      \rho(z_k,\phi_k)\rightarrow \mbox{ms} \leq
      \rho(\tilde{z},\tilde{\phi})<\infty
      \]
      From Lemma \ref{re:properties} it follows that the sequences
      $\{\phi_k\}_{k\in\N}$ and $\{z_k\}_{k\in\N}$ are uniformly bounded
      in $H^1 (\Omega)$ and $\bv(\Omega)$, respectively. Thus there exists
      subsequences, again denoted by $\{\phi_k\}_{k\in\N}$ and
      $\{z_k\}_{k\in\N}$, such that
      \[
      \phi_k \rightarrow \phi^\dag \mbox{ in } L^2(\Omega)\,,\;
      z_k\rightarrow z^\dag \mbox{ in } L^1(\Omega)\;.
      \]
      Since $\rho$ is weakly lower semi-continuous, it follows
      \[
      \mbox{ms}=\lim_{k\rightarrow\infty}\rho(z_k,\phi_k)
      \geq \rho(z^\dag,\phi^\dag)\;.
      \]
      Since $F$ is continuous on $L^1(\Omega)$,
      $F(z^\dag)=\lim_{k \to \infty} F(z_k)=y$.
      Analogously the proof of Theorem \ref{th:admissible} it follows, that
      $(z^\dag,\phi^\dag)$ is admissible and therefore
      a minimal norm solution.
\end{remunerate}\qquad
\end{proof}
\bpn Below, we summarize a stability and convergence result. The
proof uses classical techniques from the analysis of Tikhonov type
regularization methods (e.g. see
\cite{EngKunNeu89,SeiVog89,AcaVog94,EngHanNeu96,NasSch98}) and
thus is omitted:
\begin{theorem}
\label{th:con}

{\bf{Convergence:}}
      Let $\|y^\delta - y\|_{Y} \leq \delta$.
      If $\alpha = \alpha(\delta)$ satisfies
      \[
      \lim_{\delta \to 0} \alpha (\delta) = 0 \mbox{ and }
      \lim_{\delta \to 0} \frac{\delta^2}{\alpha (\delta)} = 0\;.
      \]
      Then, for a sequence $\{\delta_k\}_{k\in\N}$ converging to 0 there
      exists a sequence $\{\alpha_k := \alpha(\delta_k)\}_{k \in  \N}$
      such that $(z_{\alpha_k},\phi_{\alpha_k})$ converges in $L^1(\Omega)\times L^2(\Omega)$
      to a minimal
      norm solution.

\end{theorem}

\section{Numerical Solution}
\label{sec:2.4}
We consider a stabilized functional
\begin{equation}
\label{eq:m-reg} \Ger (\phi):= \|F(P_\ve(\phi)) - y^\delta\|_{Y}^2
+
              2\beta \alpha |P_\ve(\phi)|_{\bv} +
              \alpha \|\phi - \phi_0\|_{H^1(\Omega)}^2\;.
\end{equation}
This functional is well--posed as the following lemma shows:
\begin{lemma}
For any $\phi_0 \in H^1(\Omega)$ the functional \req{m-reg} attains a
minimizer.
\end{lemma}
\begin{proof}
The proof is similar to the proof of Theorem \ref{th:existenz-F}
by taking into account that for any sequence $\{\phi_k\}_{k \in
\N}$ converging weakly to $\phi$ in the $H^1(\Omega)$-norm, there
exists a strongly convergent subsequence in $L^2(\Omega)$.
Denoting the subsequence again by $\{\phi_k\}_{k \in \N}$ we find
\begin{remunerate}
\item \[
      \|P_\ve(\phi_k) - P_\ve(\phi)\|_{L^1(\Omega)} \leq
      \frac{1}{\ve} \mbox{meas}(\Omega)\|\phi_k - \phi\|_{\lzo} \to 0\;.
      \]
\item Therefore
      \[
      |P_\ve(\phi)|_{\bv} \leq \liminf_{k\rightarrow\infty} |P_\ve(\phi_k)|_{\bv}\;.
      \]
\end{remunerate}
Now, the assertion can be proved analogously as Theorem
\ref{th:existenz-F}.\qquad
\end{proof}
\bpn
In the following we show that for $\ve \to 0$ the minimizer of
$\Ger$ approximates a minimizer of $\Ga$, i.e., it approximates an
admissible pair.

\begin{theorem}
\label{th:just} Let $\pie$ be a minimizer of $\Ger$. Then for
$\ve_k \to 0$, there exists a convergent subsequence $(P_{\ve_k}
(\phi_{\ve_k,\alpha}),\phi_{\ve_k,\alpha}) \to
(\tilde{z},\tilde{\phi})$ in $L^1(\Omega) \times \lzo$, and the
limit minimizes $\Ga$ in the set of admissible pairs.
\end{theorem}
\begin{proof}
\begin{remunerate}
\item The infimum of $\Ga$ is attained (cf. Theorem \ref{th:admissible}),
      i.e.,
      there exists $(z_\alpha,\phi_\alpha)$ minimizing $\Ga$ over all
      admissible pairs. In particular, taking into account the definition
      of admissible pairs, there exists a sequence
      $\{\ve_k\}_{k \in \N}$ of positive numbers con\-verging to
      zero and a corresponding sequence
      $\left\{\phi_k\right\}_{k\in\N}$ in $H^1(\Omega)$
      satisfying
      \[
      \begin{aligned}
      (P_{\ve_k}(\phi_k),\phi_k) &\to (z_\alpha,\phi_\alpha)
      \mbox{ in } L^1(\Omega) \times L^2(\Omega)\,,\\
      \rho(z_\alpha,\phi_\alpha) &=
      \lim_{k \to \infty} \left\{ 2 \beta |P_{\ve_k} (\phi_k)|_\bv +
                                  \|\phi_k -\phi_0\|^2_{H^1(\Omega)}
                          \right\}\;.
      \end{aligned}
      \]
\item Let $\phi_{\ve_k}$ be a minimizer of
      ${\cal{G}}_{\ve_k,\alpha}$. The sequence
      $\{\phi_{\ve_k}\}_{k \in \N}$ is uniformly
      bounded in $H^1(\Omega)$. Thus it has a weakly convergent subsequence
      (which is again denoted by the same indices) and the weak limit is
      denoted $\tilde{\phi}$. Moreover, $\{P_{\ve_k}(\phi_{\ve_k})\}_{k\in\N}$
      is uniformly bounded in $\bv(\Omega)$. Thus, by the
      compact Sobolev embedding theorem there exists a subsequence
      $\{\phi_{\ve_k}\}_{k \in \N}$
      (again denoted with the same indices) satisfying
      \[
      \phi_{\ve_k} \to \tilde{\phi}
      \mbox{ in } L^2(\Omega)\,,
      \mbox{ and } P_{\ve_k}(\phi_{\ve_k}) \to \tilde{z}
      \mbox{ in } L^1(\Omega)\;.
      \]
      Thus $(\tilde{z},\tilde{\phi}) \in {\cal P} \times
      H^1(\Omega)$ is admissible.
\item From the definition of $\rho$ and
      the continuity of $F: L^1(\Omega) \to Y$ it
      follows that
      \[
      \begin{aligned}
      \|F(\tilde{z})-y^\delta\|^2_Y &=
      \lim_{k\rightarrow\infty} \|F(P_{\ve_k}(\phi_{\ve_k}))-y^\delta\|_Y^2\,,\\
      \rho (\tilde{z},\tilde{\phi}) &\leq
      \liminf_{k\rightarrow\infty} \left\{2\beta
      |P_{\ve_k}(\phi_{\ve_k})|_\bv+\left\|\phi_{\ve_k}-\phi_0\right\|^2_{H^1(\Omega)}\right\}
      \end{aligned}
      \]
      This shows that
      \[
      \begin{aligned}
      \Ga (\tilde{z},\tilde{\phi}) &\leq
      \liminf_{k\rightarrow\infty}
          {\cal{G}}_{\ve_k,\alpha} (\phi_{\ve_k})\\
      & \leq \liminf_{k\rightarrow\infty}
      {\cal{G}}_{\ve_k,\alpha} (\phi_k)\\
      & = \|F(z_\alpha)-y^\delta\|_Y^2 + \alpha \rho(z_\alpha,\phi_\alpha)\\
      & = \inf \Ga\;.
      \end{aligned}
      \]
      Therefore the infimum of $\Ga$ is attained at $(\tilde{z},\tilde{\phi})$.
\end{remunerate}\qquad
\end{proof}

Theorem \ref{th:just} justifies to use the functionals $\Ger$ for
approximation of the minimizer of $\Ga$.
In contrast to the minimizer of $\Ger$, which is a function in $H^1(\Omega)$,
the minimizer of $\Ga$ is an admissible pair $(z_\alpha,\phi_\alpha)$. Recall
that the function $z_\alpha$ is not uniquely defined by $\phi_\alpha$ if it
attains critical values in a neighborhood of the zero level set
(cf. Remark \ref{re:25}).

For numerical purposes it is convenient to derive the optimality
conditions of a minimizer of this functional. To this end we
consider the functional $\Ger$ with $Y=L^2(\partial \Omega)$.

Since $P'_\ve(\phi)$ is self-adjoint, we can write the formal
optimality condition for a minimizer of the functional $\Ger$
as follows:
\begin{equation}
\label{eq:formal}
\alpha (\Delta-I)(\phi - \phi_0)  = R_{\ve,\alpha,\beta} (\phi)\,,
\end{equation}
where
\[
R_{\ve,\alpha,\beta}(\phi) =
P'_\ve(\phi) F'(P_\ve(\phi))^* ( F(P_\ve(\phi)) - y^\delta )
- \beta \alpha P'_\ve(\phi) \nabla \cdot
  \left( \frac{\nabla P_\ve(\phi)}{|\nabla P_\ve (\phi)|} \right) \;.
\]

\section{Iterative Regularization and the Relation to Dynamic Level Set Methods}

For $n=1$ set ${\cal{G}}_\alpha^{(1)} (z,\phi) = {\cal{G}}_\alpha (z,\phi)$
(cf. \req{gzphi}). Iterative regularization consists in minimizing the family
of functionals
\begin{equation} \label{eq:it}
{\cal{G}}_\alpha^{(n)} (z,\phi) =
\|F(z) - y^\delta\|_Y^2 + \alpha \rho^{(n)}(z,\phi)
\end{equation}
where $\rho^{(n)}$ is the functional $\rho$ (as defined in \req{phi})
with $\phi_0$ replaced by $\phi_{n-1}$.
The minimizer of ${\cal{G}}_\alpha^{(n)} (z,\phi)$ is denoted by $\phi_n$.

Proceeding as before, we find that $\phi_n$ can be realized by solving the
formal optimality condition
\begin{equation}
\label{eq:formal-it}
\alpha (\Delta-I)(\phi - \phi_{n-1}) = R_{\ve,\alpha,\beta}(\phi)\;.
\end{equation}
Identifying $\alpha = 1/\Delta t$, $t_n = n \Delta t$, and
$\phi_n = \phi(t_n)$, $n=0,1,\ldots$ we find
\begin{equation}
\label{eq:implicit}
(\Delta-I)\left(\frac{\phi(t_n) - \phi(t_{n-1})}{\Delta t} \right) =
R_{\ve,1/\Delta t,\beta}(\phi (t_n))\;.
\end{equation}
Considering $\Delta t$ as a time discretization and using
$\beta = b_\Delta \Delta t$ we find that in a formal sense the
iterative regularized solution $\phi_n$ is
a solution of an implicit time step for the dynamic system
\begin{equation}
\label{eq:par}
(\Delta-I)\left(\frac{\partial \phi(t)}{\partial t}\right) =
R_{\ve,1/\Delta t,b_\Delta \Delta t} (\phi(t))\;.
\end{equation}
In our numerical experiments we have calculated the solution of the
dynamic system \req{par}.

\begin{figure}[h]
\centerline{\includegraphics[height=5cm,width=0.4\textwidth]{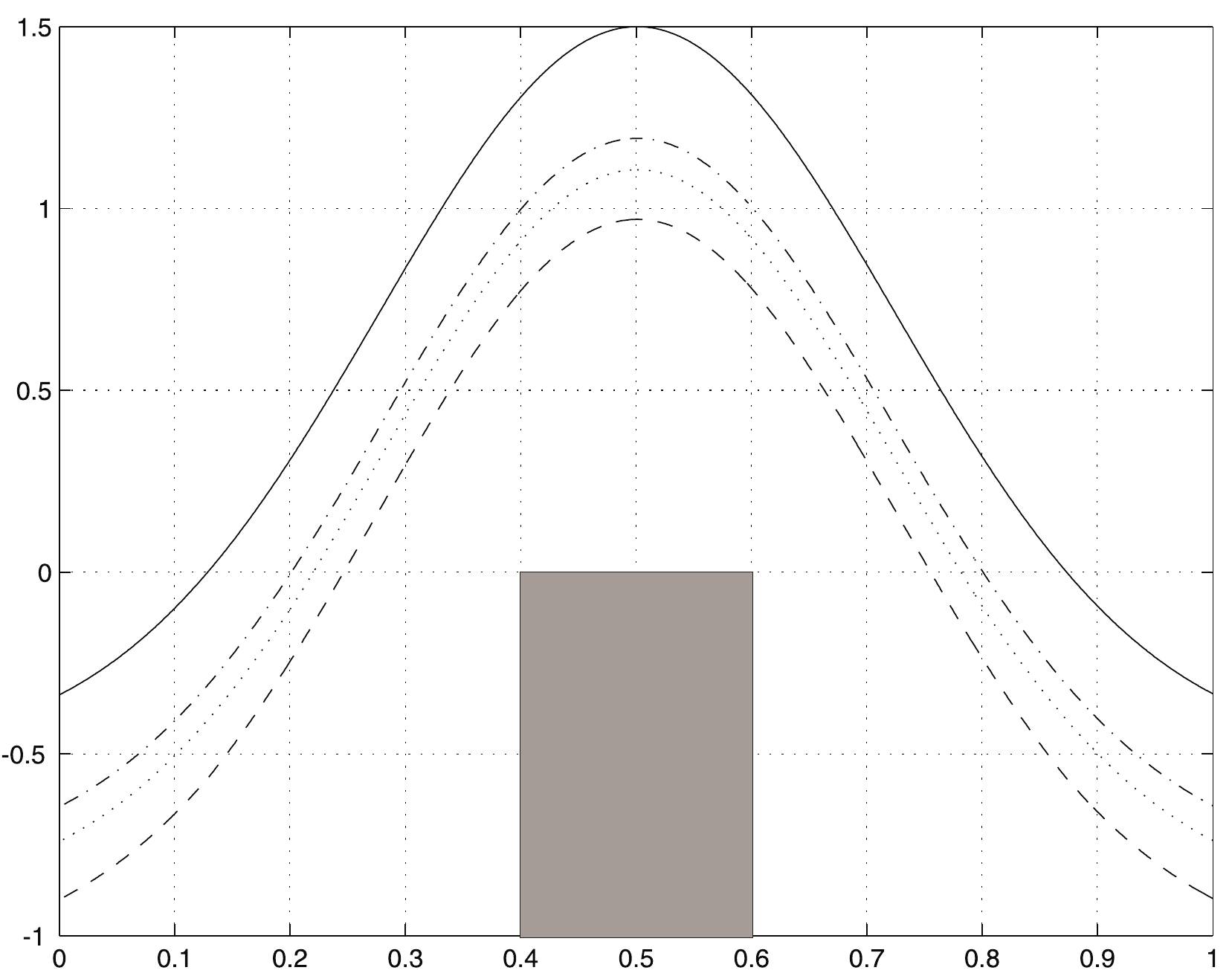}\quad
\includegraphics[height=5cm,width=0.4\textwidth]{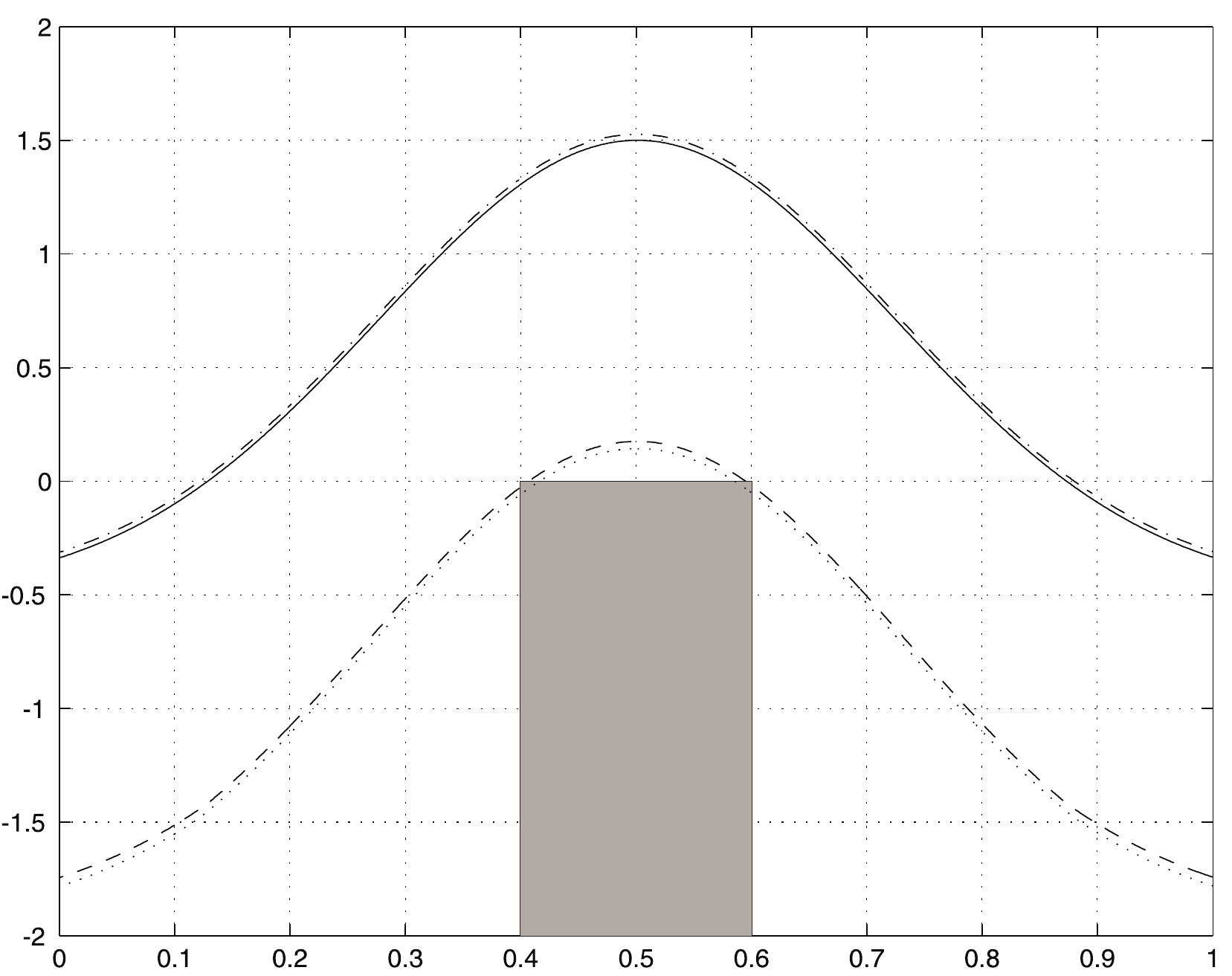}}

\caption{\label{shift} The functions $\phi^{(0)}$ (solid line),
$\phi^{(1)}$ (dashed line), $\phi^{(2)}$ (dash-dot line) and
$\phi^{(3)}$ (dotted line). To recover is the interval $[0.4,
0.6]$, which is displayed by the grey rectangle. The first
iteration is the best. In the right picture $\alpha$ is smaller
than in the left picture.}
\end{figure}

\begin{figure}[h]
\centerline{\includegraphics[height=5cm,width=0.8\textwidth]{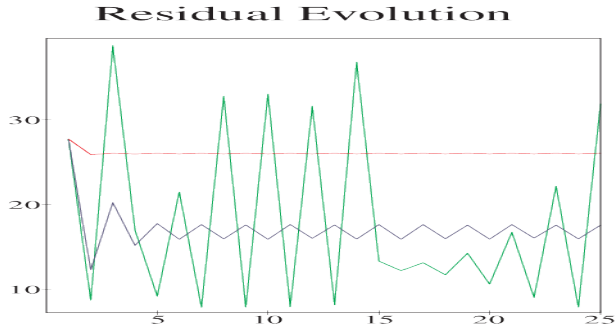}}
\caption{\label{fig:neu}
Decay of the residual $\|F(P_\ve(\phi^{(k)})) - y^\delta\|_Y^2$ in dependence
of the number of iterations (residual evaluated for the first experiment --
noise free data in Section \ref{sec:2.5}). After the first iteration the
fixed-point iteration stagnates.}
\end{figure}

\begin{figure}[h]
\centerline{\includegraphics[height=5cm,width=0.8\textwidth]{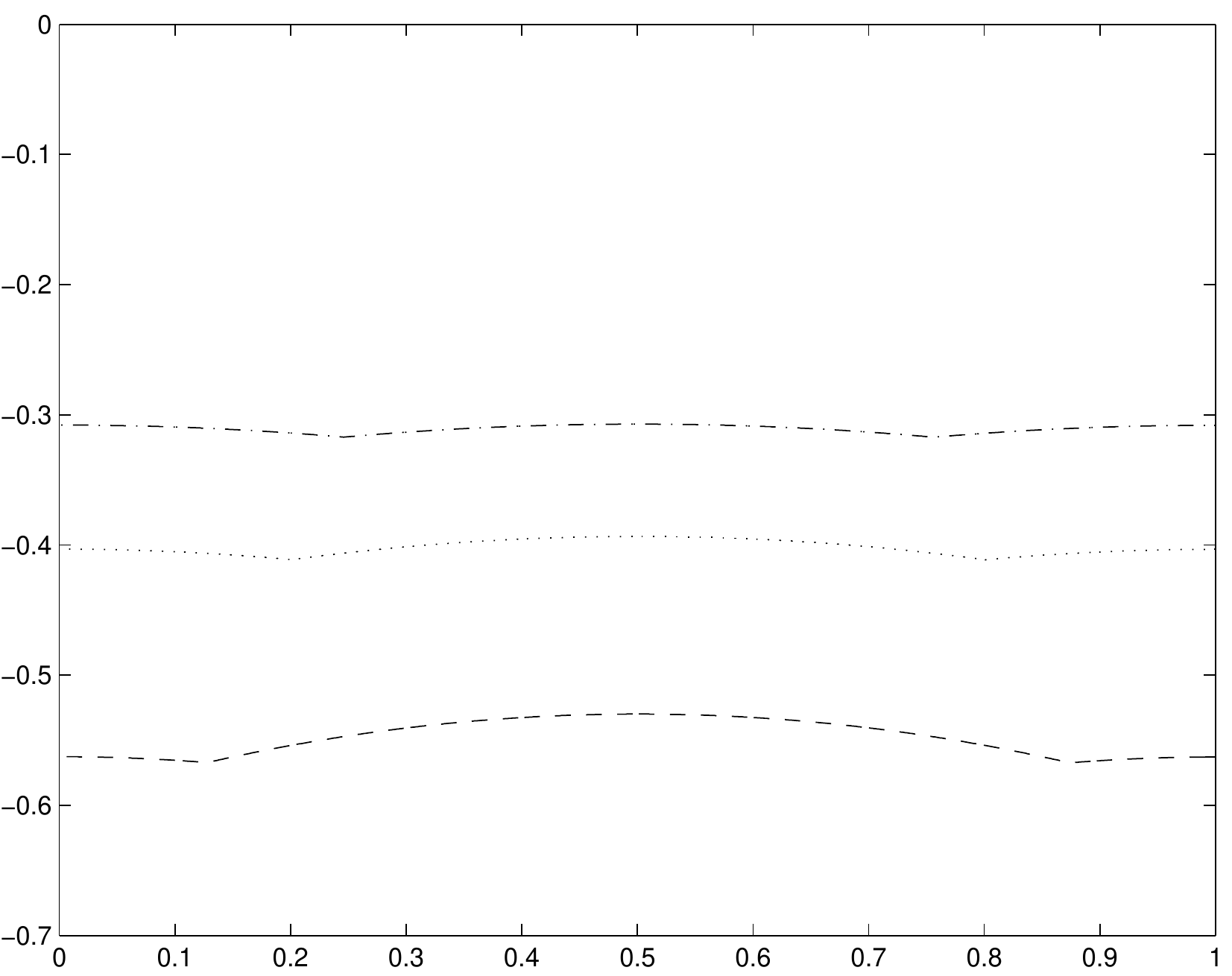}}
\caption{\label{osziallation} The differences between $\phi^{(0)}$
and the functions $\phi^{(1)}$ (dashed line), $\phi^{(2)}$
(dash-dot line) and $\phi^{(3)}$ (dotted line) from the left
picture of figure \ref{shift}.}
\end{figure}

For each time step it is required to solve equation \req{implicit}.
$\phi(t_n)$ in \req{implicit} can be solved with a fixed point iteration:
setting $\phi(t_{n-1}) = \phi^{(0)}$, we get
$\phi(t_n) = \lim_{k \to \infty} \phi^{(k)}$
\begin{equation}
\label{eq:stag} (\Delta-I)\left(\frac{\phi^{(k+1)} -
\phi^{(0)}}{\Delta t} \right) = R_{\ve,1/\Delta t,b_\Delta \Delta
t} (\phi^{(k)})\;.
\end{equation}
In our numerical experiments we observed that the iteration does
not significantly improve after the first iteration (cf. Figure
\ref{shift}).  This behavior becomes transparent by noting the
$H^1$-seminorm typically dominates the $L^2$-norm in the quadratic
regularization term. The $H^1$-seminorm difference of the
regularized solution and $\phi^{(0)}$ is small if it is just
shifted up or down. In numerical experiments it is observed that
the first iteration almost corresponds to a horizontal shift of
$\phi^{(0)}$ such that the residual functional is minimized (cf.
Figure \ref{fig:neu}) and also the further iterations are again
nearly horizontally shifted versions of $\phi^{(0)}$ (cf. Figure
\ref{osziallation}).

In almost all test examples the residual $\|F(P_\ve(\phi^{(k)})) -
y^\delta\|^2$ is oscillating in dependence of $k$ (cf. Figure
\ref{fig:neu}) and smallest for $k=1$.

The above consideration justify to restrict attention to the approximate
solution of the dynamic system \req{formal-it} where in each time step only
one iteration step of \req{stag} is used, i.e., we use an explicit
Euler method for solving the evolution process.
In this case numerical instabilities may occur by dividing by small
absolute values of the gradient in the differential
$\nabla \cdot
  \left( \frac{\nabla P_\ve(\phi)}{|\nabla P_\ve (\phi)|} \right)$.
Thus, for numerical purpose it is convenient to introduce a small positive
number $h$ and replace the differential by
\[
\nabla \cdot
  \left( \frac{\nabla P_\ve(\phi)}{\sqrt{|\nabla P_\ve (\phi)|^2+h^2}} \right)\;.\]
Usually semi--implicit iteration schemes require a less
restrictive time marching (this approach is commonly referred as
Dziuk's method). The implementation would require to solve
\begin{equation}
\label{eq:stag*}
\begin{aligned}
(\Delta-I)\left(\frac{\phi^{(k+1)} - \phi^{(0)}}{\Delta t} \right)
& = P'_\ve(\phi^{(k)})
   F'(P_\ve(\phi^{(k)}))^* ( F(P_\ve(\phi^{(k)})) - y^\delta ) \\
& \quad
- b_\Delta  P'_\ve(\phi^{(k)}) \nabla \cdot
  \left( \frac{\nabla P_\ve(\phi^{(k+1)})}{\sqrt{|\nabla P_\ve (\phi^{(k)})|^2+h^2}}
   \right) \;.
\end{aligned}
\end{equation}
In implementation of this approach the difficulty arises that the
function in front of
$\nabla \cdot
  \big( \nabla P_\ve(\phi^{(k+1)}) /
         \sqrt{|\nabla P_\ve (\phi^{(k)})|^2+h^2} \big)$
vanishes outside of a neighborhood of the zero level set,
which makes it almost impossible to implement this scheme efficiently.

\section{Numerical Experiments}
\label{sec:2.5} In this section we shall consider an inverse
potential problem of recovering the shape of a domain $D$ using
the knowledge of its (constant) density and the measurements of
the Cauchy data of the corresponding potential on the boundary of
a fixed Lipschitz domain $\Omega \subset \R^2$, which contains
$\overline D$. This is the same problem as considered by Hettlich
and Rundell \cite{HetRun96}, which used iterative methods for
recovering a single star-shaped object.

To achieve an analogous problem, a certain definition of the
operator $F$ is necessary:
\[ \begin{aligned}
    F: \lzo & \rightarrow & L^2(\partial \Omega)\\
    \chi_D & \rightarrow & F(\chi_D)
    \end{aligned}
\]
This is possible, because we consider only characteristic
functions $\chi_D$. The $\lzo$-norm is then equivalent to the
$L^1(\Omega)$-norm of $\chi_D$. Therefore the necessary properties
are retained.

 The problem introduced above can mathematically be described
as follows:
\begin{equation} \label{eq:direct-probl}
\Delta u = \chi_D\, ,\ {\rm in}\ \Omega\,;\ \ u|_{\partial\Omega} = 0\, ,
\end{equation}
where $\chi_D$ is the characteristic function of the domain $D
\subset \Omega$, which has to be reconstructed. Since $\chi_D \in
\lzo$, the Dirichlet boundary value problem in \req{direct-probl}
has a unique solution, the potential $u \in H^2(\Omega) \cap
H^1_0(\Omega)$. Here $H^1_0(\Omega)$ is defined as the closure
with respect to $H^1(\Omega)$ of functions in $C^\infty(\Omega)$
with compact support in $\Omega$.

The inverse problem we are concerned with, consists in determining the
shape of $D$ from measurements of the Neumann trace of $u$ at
$\partial\Omega$, i.e. from $[\partial u /\partial\nu]_{\partial\Omega}$,
where $\nu$ represents the outer normal vector to $\partial\Omega$.

Notice that this problem can be considered in the framework of an
inverse problem for the {\em Dirichlet to Neumann map}. For given
$h \in \lzo$, the Dirichlet to Neumann operator maps a Dirichlet
boundary data onto the Neumann trace of the potential, i.e.,
$\Lambda: H^{1/2}(\partial\Omega) \to H^{-1/2}(\partial\Omega)$,
$\Lambda(\varphi) := [\partial \tilde{u} /
\partial\nu]_{\partial\Omega}$, where $\tilde{u}$ solves
\[
\Delta \tilde{u} = h\, ,\ {\rm in}\ \Omega\,;\ \
\tilde{u}|_{\partial\Omega} = \varphi\, .
\]
The inverse problem for the $\Lambda$ operator consists in
determining the unknown parameter (i.e., the function $h$) from
different pairs of Dirichlet, Neumann boundary data. The general
case with $h \in \lzo$ has already been considered by many
authors, among them we mention \cite{CabBel71,Rin95}, which
introduced numerical methods based on Tikhonov regularization, and
\cite{HetRun96} with iterative regularization methods.

Hettlich and Rundell \cite{HetRun96} observe that, in the
particular case $h = \chi_D$, one pair of measurement data of
Dirichlet--Neumann data furnishes as many information as the full
Dirichlet--Neumann operator, i.e., it is sufficient to consider
only one pair of Cauchy data for the inverse problem.
Therefore, no further information on $D$ can be gained by using various
pairs of Dirichlet--Neumann data, since we can always reduce the
reconstruction  problem to the homogeneous Dirichlet case.

For the particular case $h = \chi_D$, it has been observed by
Hettlich and Rundell \cite{HetRun96} that the Cauchy data may not
furnish enough information to reconstruct the boundary of $D$,
e.g., if $D$ is not simply connected. On the other hand, Isakov
observed in \cite{Isa90} that star like domains $D$ are uniquely
determined by their potentials.

The inverse potential problem is discussed within the general
framework introduced in Section~\ref{sec:introd}. In particular,
we allow domains, that consists of a number of connected
inclusions. For this general class we have not unique
identifiability and we restrict attention to ``minimum-norm
solutions''. Recall that in this case a minimum-norm-solution is a
level set function $\phi$, where $P (\phi)$ determines the
inclusion. A minimum norm solution satisfies that it minimizes
the functional $\rho(z,\phi)$ in the class of level set functions
such that the according Neumann boundary values $\frac{\partial
u}{\partial \nu}$ fit the data $y^\delta$.

\subsection{The level set regularization algorithm}

In the following we describe the level set regularization
algorithm. This method compares to the Landweber iteration as
proposed by Hettlich and Rundell \cite{HetRun96}. In our context
the operator $F'$ can be considered as an approximation of the
\emph{domain derivative operator}  for multiple connected domains
(cf. Figure \ref{fig:algor}).

The complexity of our algorithm is as follows: at each
iteration of the level set method, three elliptic boundary value
problems are solved (two of Dirichlet type and one of Neumann
type).

In Table~\ref{fig:algor} the iteration procedure for the solution
of the formal optimality condition \req{formal} is outlined. The
algorithm can be implemented using finite element codes (as we
did) or finite difference methods for the solution of partial
differential equations.
\begin{figure}
\hfil
\fbox{
\begin{minipage}{\textwidth}
\begin{tt}\begin{remunerate}
\item[{\bf 1.}] \ Evaluate the residual\, $r_k := F(P_\ve
(\phi_k))
- y^\delta = \frac{\partial u_k}{\partial \nu} - y^\delta$,\\
where $u_k$ solves
$$ \Delta u_k = P_\ve (\phi_k)\, ,\ {\tt in}\ \Omega\,;\qquad
   u_k|_{\partial\Omega} = 0\, . $$
\item[{\bf 2.}] \ Evaluate\, $v_k := F'(P_\ve (\phi_k))^*(r_k) \in
\lzo$, solving
$$ \Delta v_k = 0\, ,\ {\tt in}\ \Omega\,;\quad
    v_k|_{\partial\Omega} = r_k\, . $$
\item[{\bf 3.}] \ Evaluate $w_k \in H^1(\Omega)$, satisfying
\[\begin{aligned}
(I - \Delta) w_k &= - P'_\ve(\phi_k) \, v_k +
   \beta \alpha P'_\ve(\phi_k) \nabla \cdot
   \left( \frac{\nabla P_\ve(\phi_k)}{|\nabla P_\ve (\phi_k)|} \right)\,,
   {\tt in}\ \Omega;\\
\frac{\partial w_k}{\partial \nu}|_{\partial\Omega}&=0\;.
  \end{aligned}
\]
\item[{\bf 4.}] \ Update the level set function\, $\phi_{k+1} =
\phi_k + \frac{1}{\alpha} \; w_k$.
\end{remunerate} \end{tt}
\end{minipage} }
\hfil \caption{Implementation of a single iteration step for minimizing the
level set regularization. \label{fig:algor} }
\end{figure}

\subsection{Reconstruction of a density function with non simply
connected support} \label{ssec:exper1}

\begin{figure}[h]
\centerline{
\includegraphics[width=0.45\textwidth]{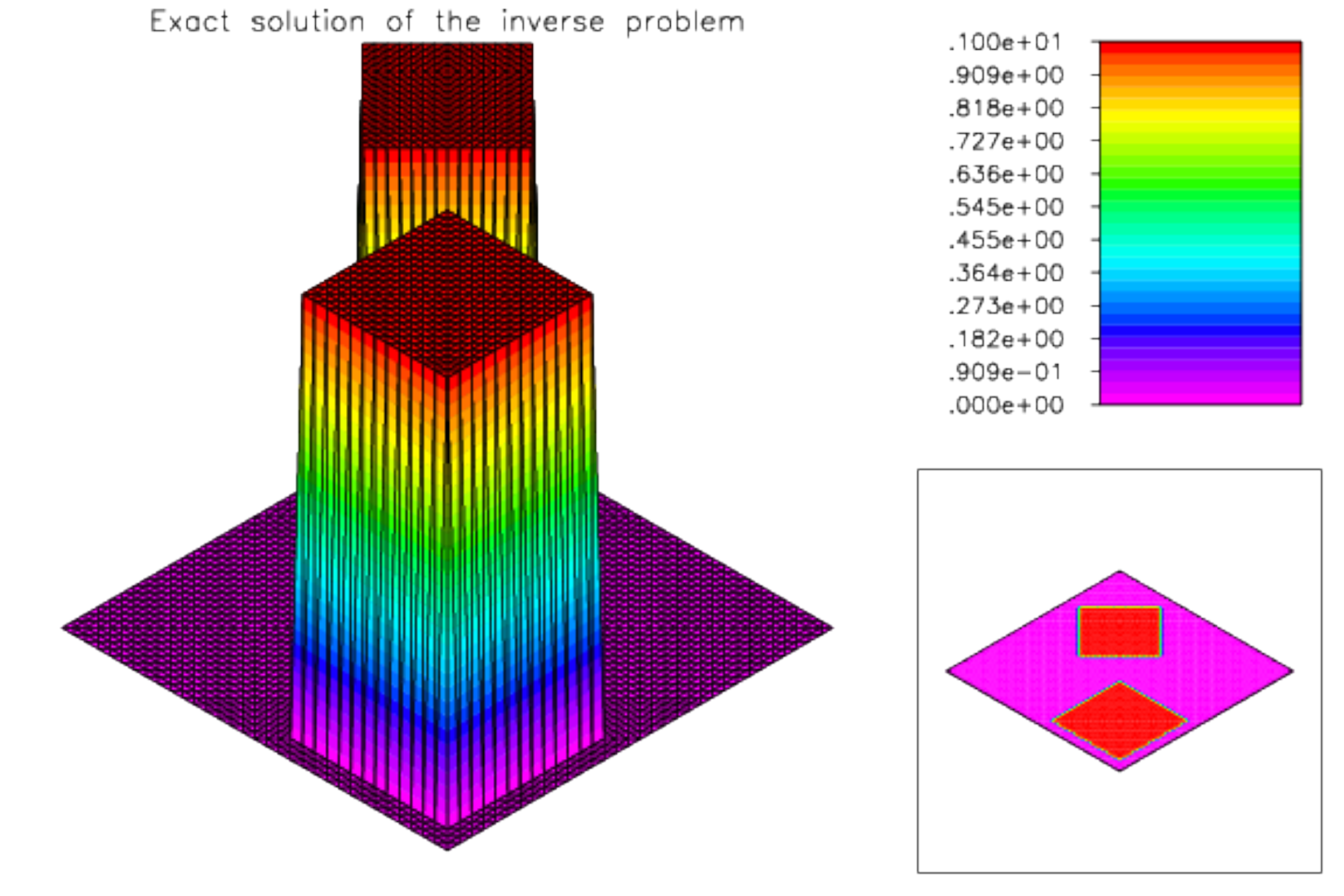} \hskip1cm
\includegraphics[width=0.45\textwidth]{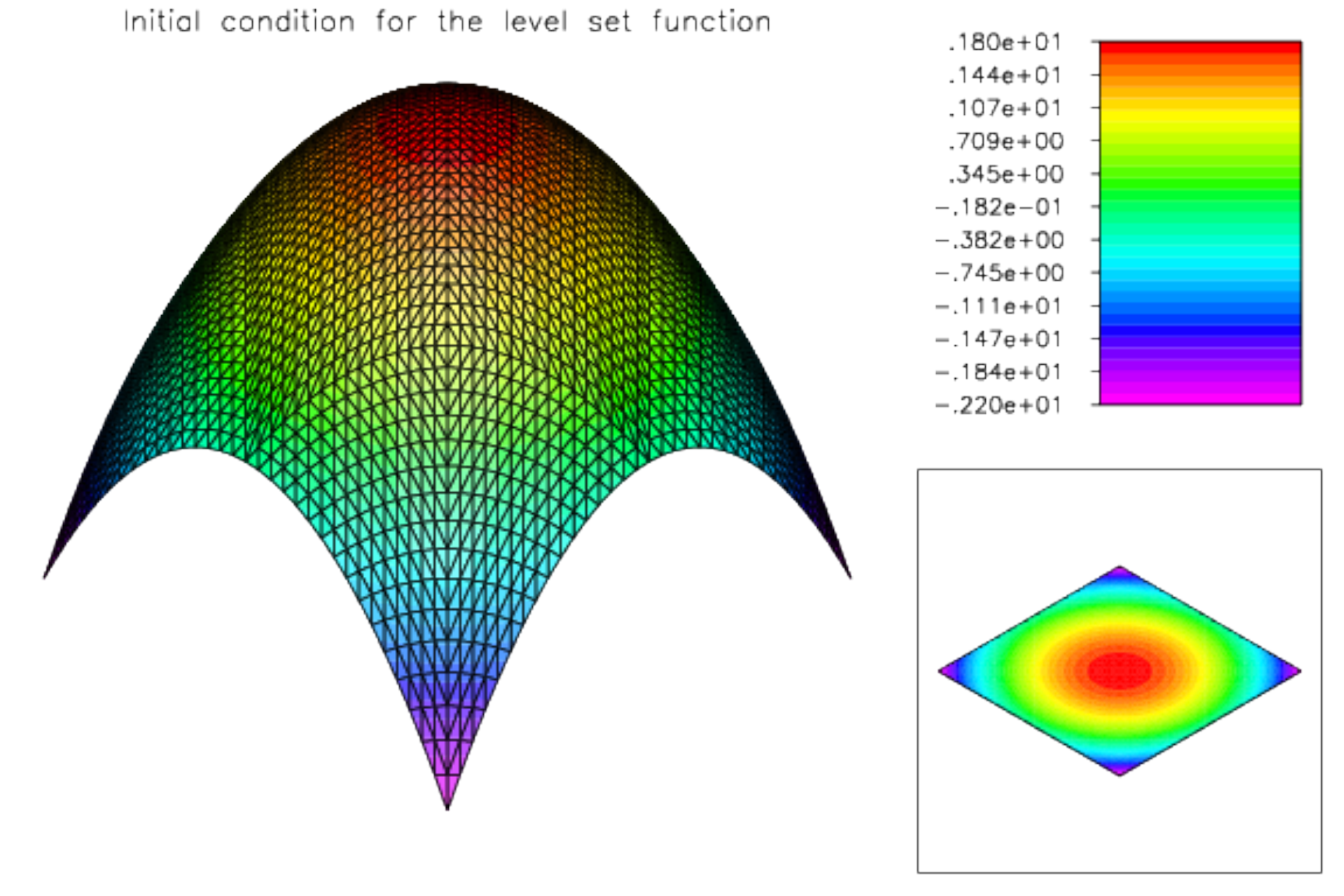}}
\caption{The picture on the left hand side shows the coefficient to
be reconstructed. On the other picture, the initial condition for
the level set regularization method. \label{fig:exsol-incond} }
\end{figure}

In this first experiment we consider the inverse problem of
reconstructing the right hand side $\chi_D$ in \req{direct-probl}
from the knowledge of a single pair of boundary data
$(u, \Lambda u)=(0,y^\delta)$ at $\partial\Omega$.
In the examples considered below we always use the squared domain
$\Omega=(0,1)^2 \subset \R^2$. $\chi_D \in \lzo$ is the characteristic
function as represented in Figure~\ref{fig:exsol-incond}.

The overdetermined boundary measurement data $y^\delta$
for solving the inverse pro\-blem, is obtained
by solving the elliptic boundary value problem 
in \req{direct-probl}. Notice that $\chi_D$ corresponds to the
characteristic function of a not-connected proper subset of
$\Omega$.
The initial condition for the level set function is shown in Figure~%
\ref{fig:exsol-incond}.
In order to avoid inverse crimes, the direct problem
\req{direct-probl} is solved on an adaptively refined grid with
8.807 nodes (three levels of adaptive refinement). Alternatively,
in the numerical implementation of the level set method, all
boundary value problems are solved at an uniformly refined grid
with 2.113 nodes.

When the data is given exactly, we tested the iterative level set
regularization without the additional regularization term
$|P_\ve(\phi_k)|_\bv$, i.e. $\beta = 0$.

In all computed experiments we use the operator $P_\ve$ defined in
Section~\ref{sec:lsr} with $\ve = 1/8$. This seams to be
compatible with the size of our mesh, since the diameter of the
triangles in the uniform grid (used in the finite element method)
is approximately $\sqrt{2}/32$.

\begin{figure}[h]
\centerline{
\includegraphics[width=0.24\textwidth]{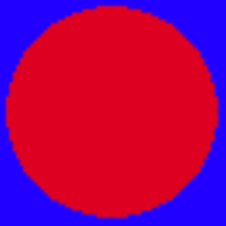}
\includegraphics[width=0.24\textwidth]{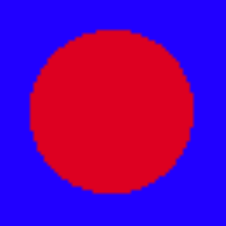}
\includegraphics[width=0.24\textwidth]{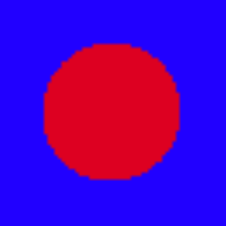}
\includegraphics[width=0.24\textwidth]{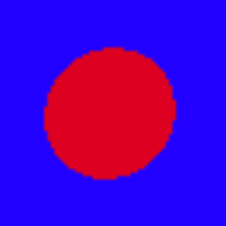} } \vskip2pt
\centerline{
\includegraphics[width=0.24\textwidth]{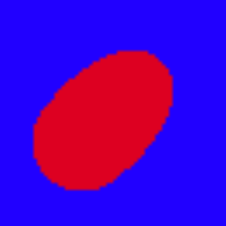}
\includegraphics[width=0.24\textwidth]{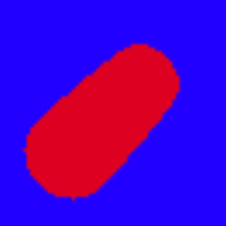}
\includegraphics[width=0.24\textwidth]{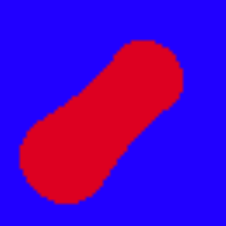}
\includegraphics[width=0.24\textwidth]{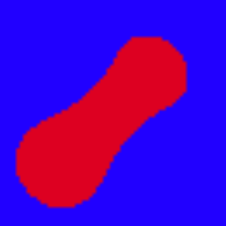} } \vskip2pt
\centerline{
\includegraphics[width=0.24\textwidth]{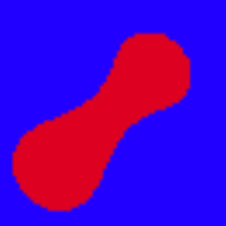}
\includegraphics[width=0.24\textwidth]{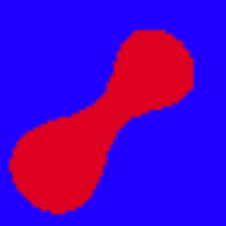}
\includegraphics[width=0.24\textwidth]{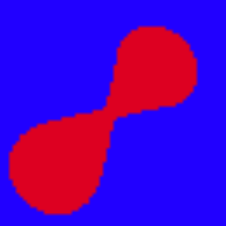}
\includegraphics[width=0.24\textwidth]{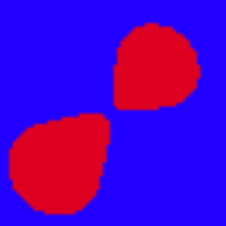} } \vskip2pt
\centerline{
\includegraphics[width=0.24\textwidth]{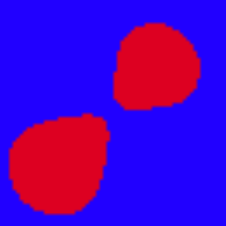}
\includegraphics[width=0.24\textwidth]{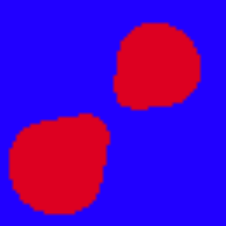}
\includegraphics[width=0.24\textwidth]{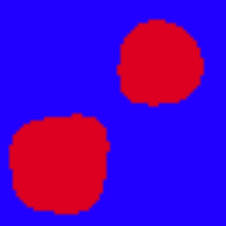}
\includegraphics[width=0.24\textwidth]{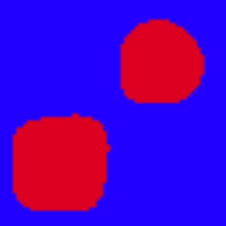} }
\caption{Level set evolution for exact data. Plots after $0,1,2,10$,
$100,200,300,400$, \ $500,600,700,800$, \ $900,1000,2000,3000$ iterative
steps. \label{fig:ls-evol-1} }
\end{figure}

In Figure~\ref{fig:ls-evol-1} we present the evolution of the
level set function for given exact data for the first 3000
iterative steps. As one can see in this figure, the original level
set splits into two convex components after approximately 800
iterations. After 1000 iterations the level set function still
changes, but very slowly. We performed similar tests for different
initial conditions and observed that, after 1000 iterations, the
corresponding pictures look very much alike.

For the second part of this experiment, the density function to be
reconstructed is still the one shown in
Figure~\ref{fig:exsol-incond}. This time, however, we add randomly
generated noise to the data $y^\delta$
used in the first part of the experiment.

\begin{figure}[h]
\centerline{\ \ \ \
\includegraphics[height=4.6cm,width=0.52\textwidth]{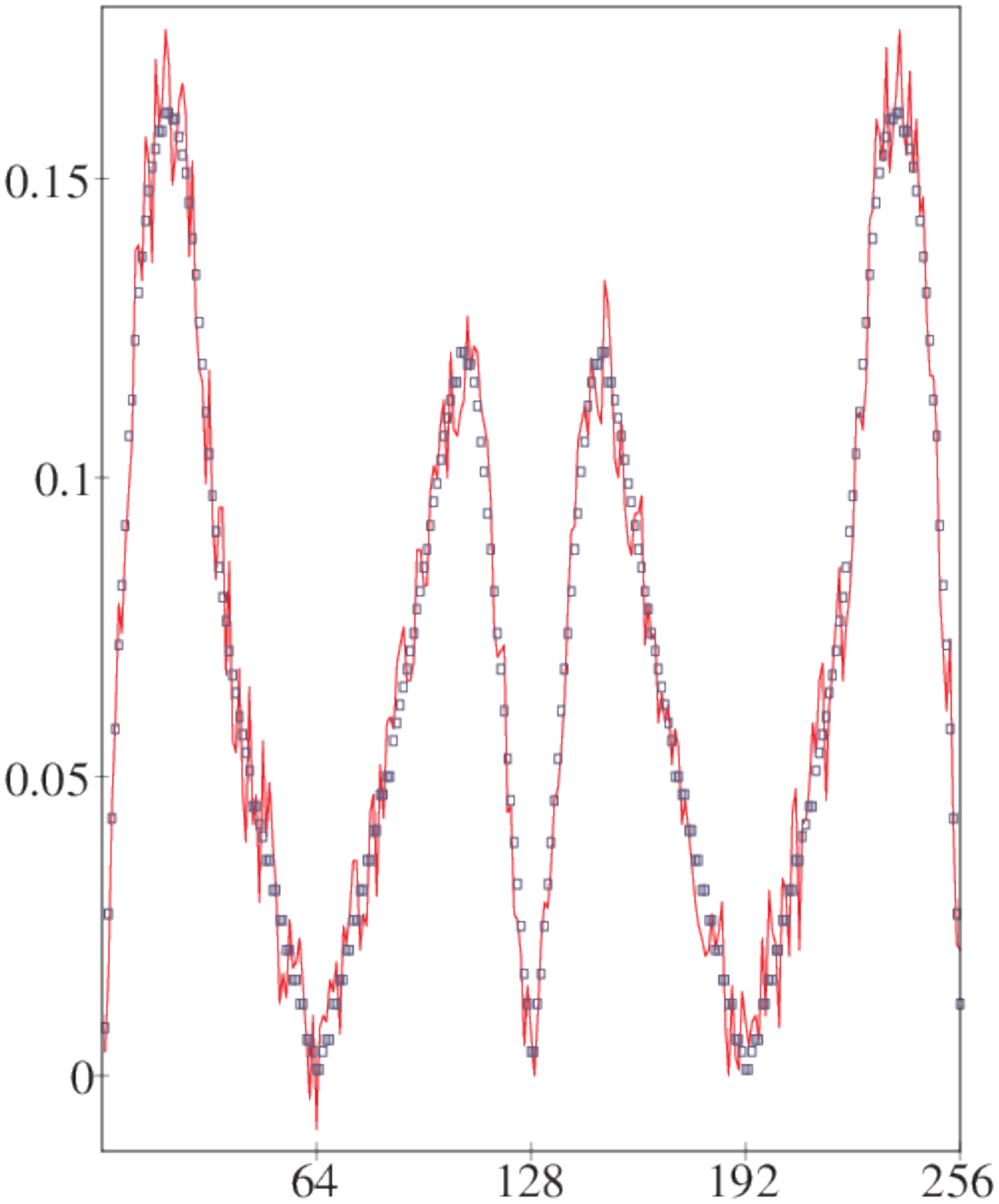}
\hskip0.4cm
\includegraphics[height=4.6cm,width=0.52\textwidth]{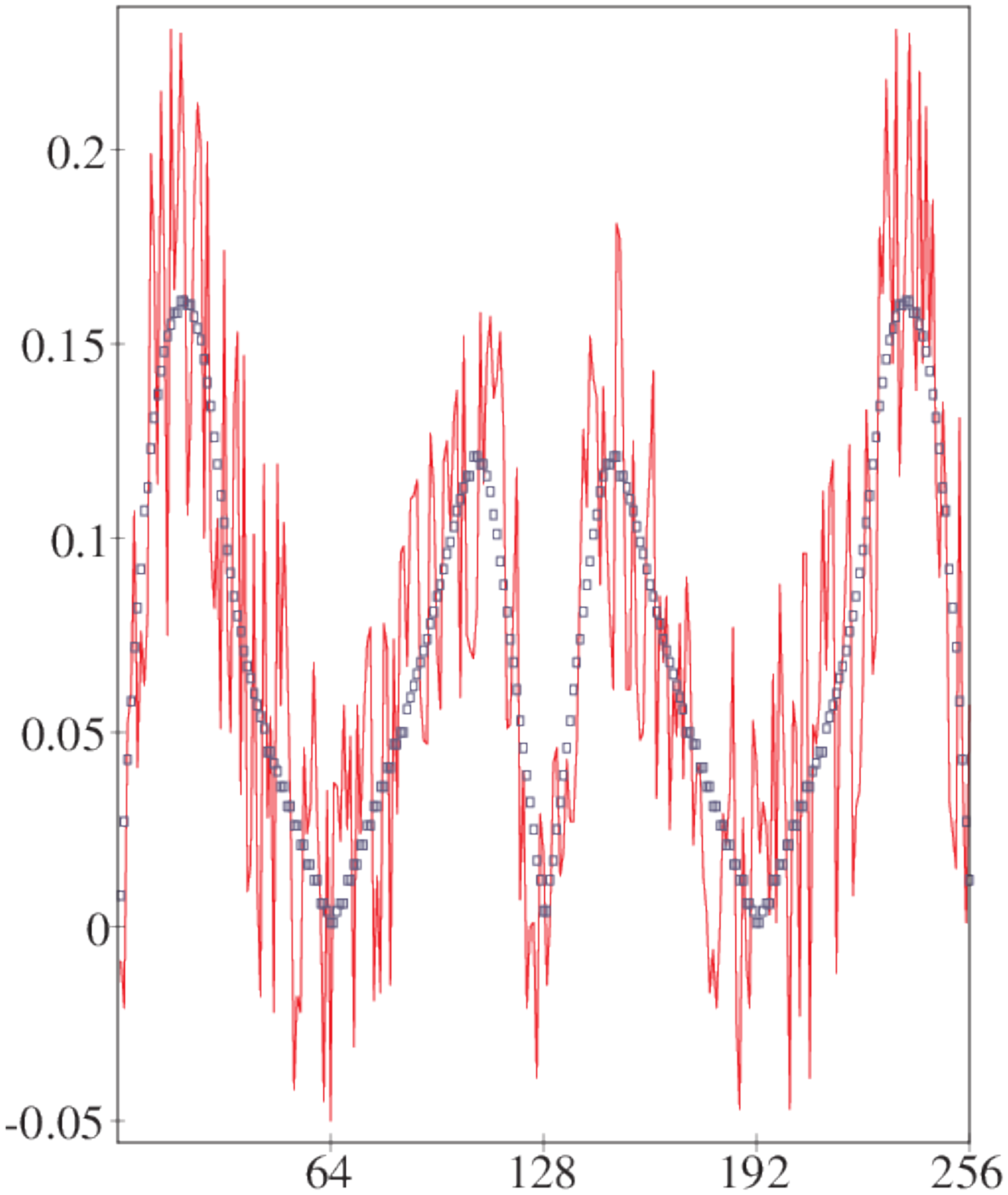}}
\vskip-.4cm \caption{The dotted (blue) line represents the exact
data $y^\delta$; 
 the solid (red) line represents the perturbed
data. The noise level corresponds to 10\% at the left hand side
and 50\% at the right hand side. \label{fig:noise} }
\end{figure}

The exact boundary data $y^\delta$
is shown in Figure~\ref{fig:noise} as the dotted (blue) line. We
consider actually two distinct sets of perturbed data: For the
first experiment we add to the exact data a white noise of $10\%$
(in the $l_\infty$-norm); For the second experiment we use a noise
level of $50\%$. Both sets of inaccurate data are plotted in
Figure~\ref{fig:noise} and correspond to the solid (red) line.

As in the noise free experiment, the same care was taken to avoid
inverse crimes. The choice of the parameter $\ve$ (operator
$P_\ve$) follows also the same criteria as before. However, since
we are now dealing with noisy data, we have to develop a strategy
for the choice of the regularization parameter $\beta$. For this
proposal we opted for the fit-to-data strategy, i.e. $\beta
\alpha$ is chosen such that the regularization term (see
Figure~\ref{fig:algor}) has the same order as the noise level.

\begin{figure}[h]
\centerline{
\includegraphics[width=0.24\textwidth]{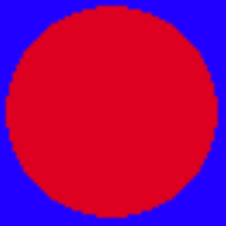}
\includegraphics[width=0.24\textwidth]{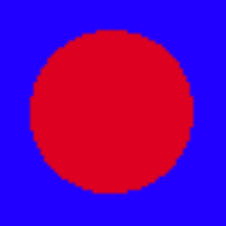}
\includegraphics[width=0.24\textwidth]{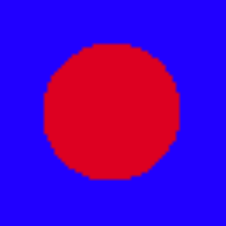}
\includegraphics[width=0.24\textwidth]{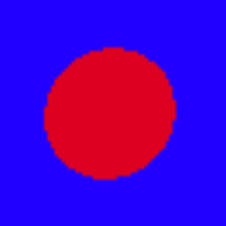} } \vskip2pt
\centerline{
\includegraphics[width=0.24\textwidth]{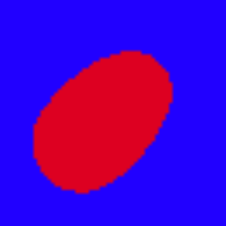}
\includegraphics[width=0.24\textwidth]{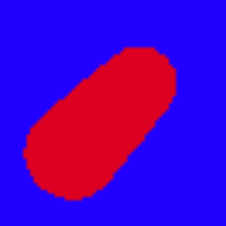}
\includegraphics[width=0.24\textwidth]{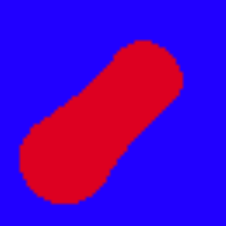}
\includegraphics[width=0.24\textwidth]{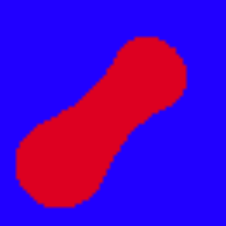} } \vskip2pt
\centerline{
\includegraphics[width=0.24\textwidth]{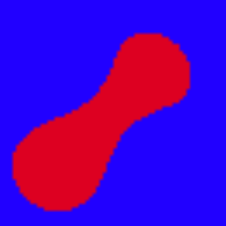}
\includegraphics[width=0.24\textwidth]{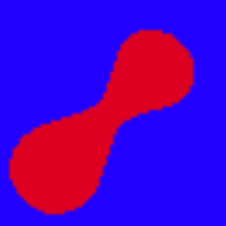}
\includegraphics[width=0.24\textwidth]{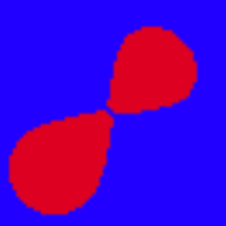}
\includegraphics[width=0.24\textwidth]{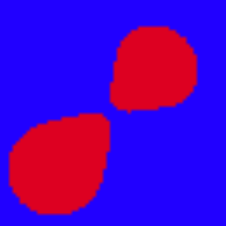} } \vskip2pt
\centerline{
\includegraphics[width=0.24\textwidth]{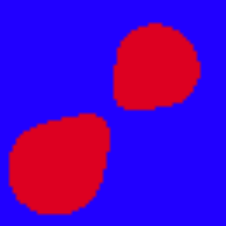}
\includegraphics[width=0.24\textwidth]{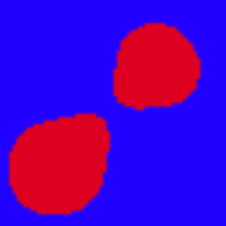}
\includegraphics[width=0.24\textwidth]{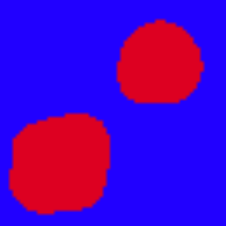}
\includegraphics[width=0.24\textwidth]{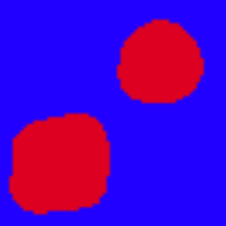} }
\caption{Level set evolution for inaccurate data; noise level of $10\%$.
Plots after $0,1,2,10$, $100,200,300,400$, $500,600,700,800$,
$900,1000,2000,3000$ iterative steps. \label{fig:ls-evol-2} }
\end{figure}

\begin{figure}[h]
\centerline{
\includegraphics[width=0.24\textwidth]{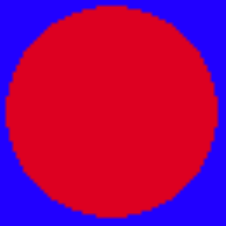}
\includegraphics[width=0.24\textwidth]{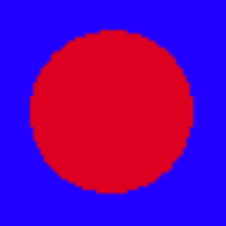}
\includegraphics[width=0.24\textwidth]{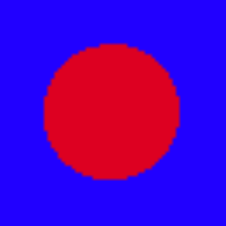}
\includegraphics[width=0.24\textwidth]{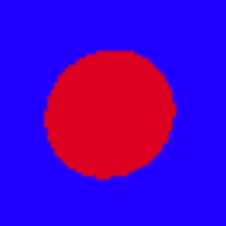} } \vskip2pt
\centerline{
\includegraphics[width=0.24\textwidth]{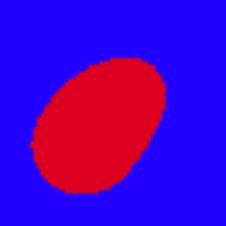}
\includegraphics[width=0.24\textwidth]{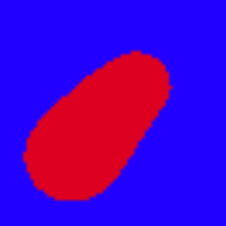}
\includegraphics[width=0.24\textwidth]{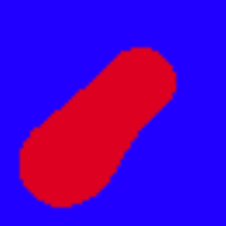}
\includegraphics[width=0.24\textwidth]{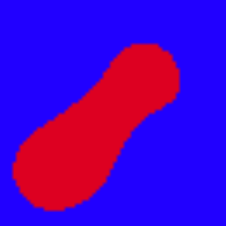} } \vskip2pt
\centerline{
\includegraphics[width=0.24\textwidth]{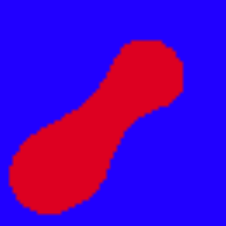}
\includegraphics[width=0.24\textwidth]{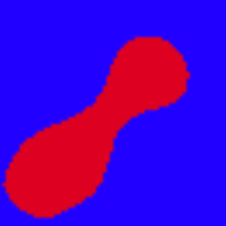}
\includegraphics[width=0.24\textwidth]{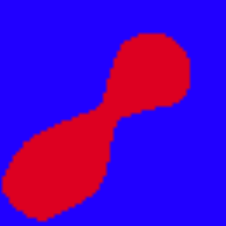}
\includegraphics[width=0.24\textwidth]{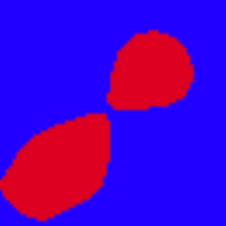} } \vskip2pt
\centerline{
\includegraphics[width=0.24\textwidth]{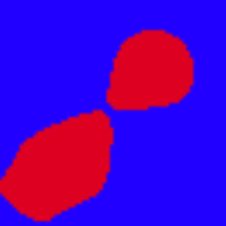}
\includegraphics[width=0.24\textwidth]{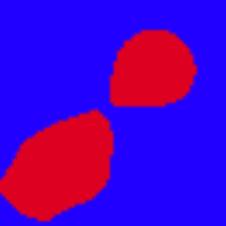}
\includegraphics[width=0.24\textwidth]{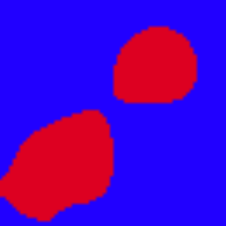}
\includegraphics[width=0.24\textwidth]{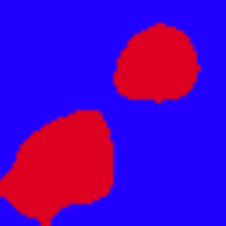} }
\caption{Level set evolution for inaccurate data; noise level of $50\%$.
Plots after $0,1,2,10$, $100,200,300,400$, $500,600,700,800$,
$900,1000,1300,1600$ iterative steps. \label{fig:ls-evol-3} }
\end{figure}

The corresponding results generated by the level set method where
surprisingly stable, as one can observe in
Figures~\ref{fig:ls-evol-2} and~\ref{fig:ls-evol-3}. In the first
case (noise level of $10\%$) the results are comparable with the
previous experiment, where exact data was available. In the second
case (noise level of $50\%$) we are not able to precisely recover
the shape of the set $D$, corresponding to the characteristic
function shown in Figure~\ref{fig:exsol-incond}. However, we are
still able to identify the number of connected components of $D$,
as well as their relative positions inside the domain $\Omega$.

\subsection{Reconstruction of a density function with non convex support}

\begin{figure}[h]
\centerline{
\includegraphics[width=0.45\textwidth]{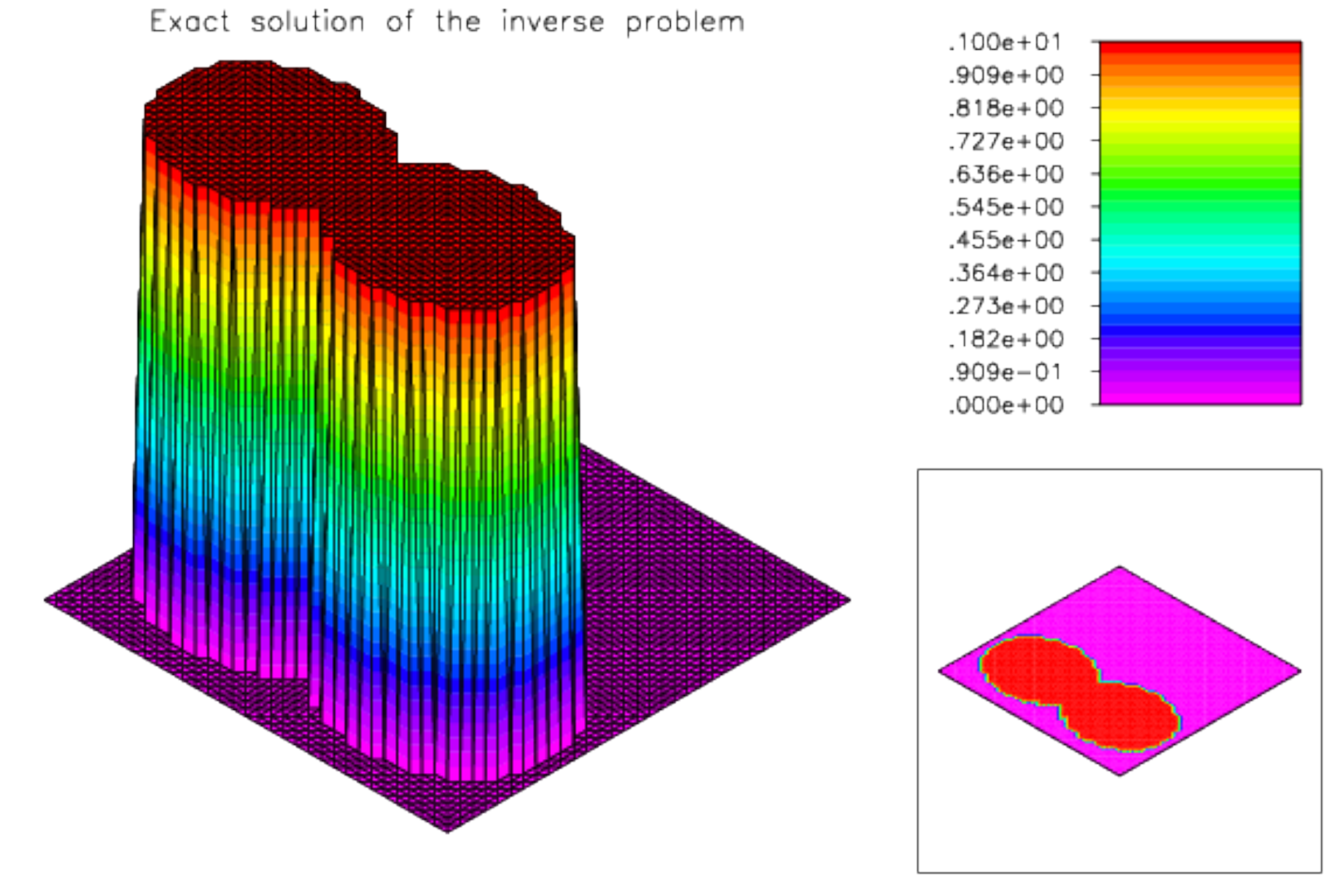} \hskip1cm
\includegraphics[width=0.45\textwidth]{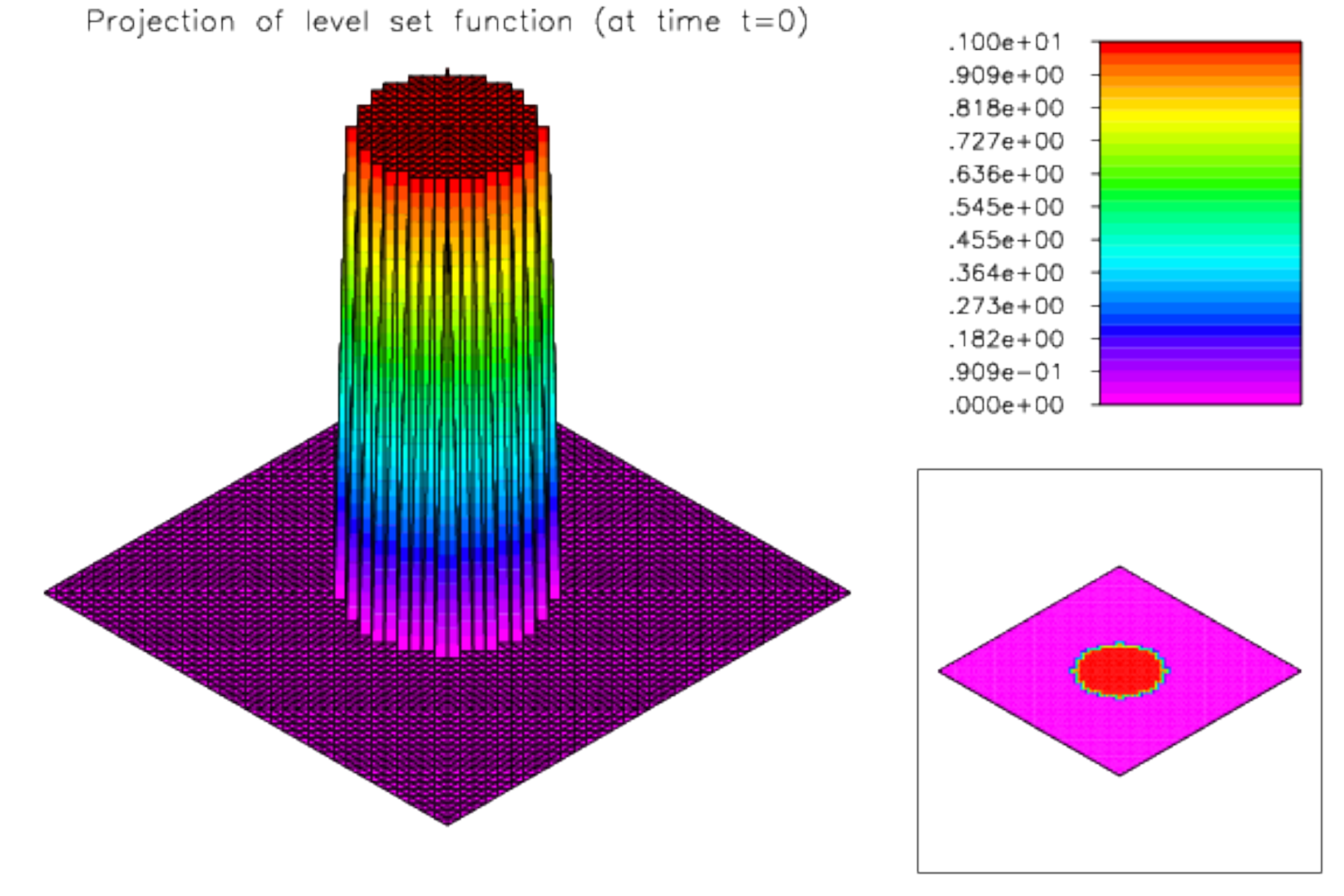}}
\caption{The picture on the left hand side shows the coefficient
to be reconstructed. On the other picture, the (projection of the)
initial condition for the level set regularization method.
\label{fig:exsol-incond2}}
\end{figure}

In this second experiment we consider the problem of reconstructing the
density function shown in Figure~\ref{fig:exsol-incond2}. The main goal
now is to investigate the difficulty of the level set method in recovering
non convex domains. The domain $\Omega$ is the same used in
Subsection~\ref{ssec:exper1} and again we aim to reconstruct the density
function in (\ref{eq:direct-probl}) from boundary measurements.

\begin{figure}[ht]
\centerline{
\includegraphics[width=0.24\textwidth]{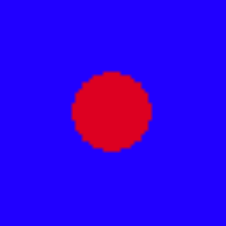}
\includegraphics[width=0.24\textwidth]{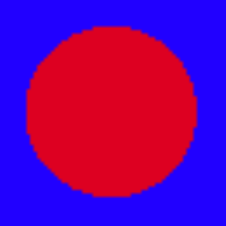}
\includegraphics[width=0.24\textwidth]{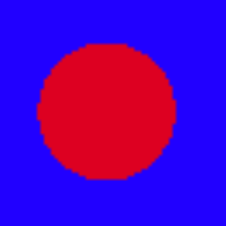}
\includegraphics[width=0.24\textwidth]{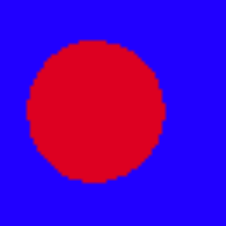} } \vskip2pt

\centerline{
\includegraphics[width=0.24\textwidth]{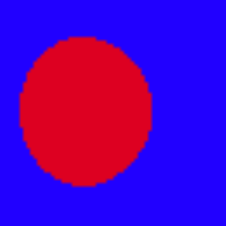}
\includegraphics[width=0.24\textwidth]{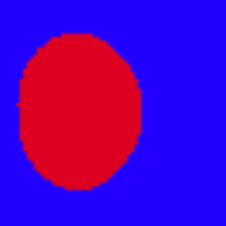}
\includegraphics[width=0.24\textwidth]{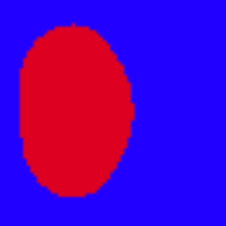}
\includegraphics[width=0.24\textwidth]{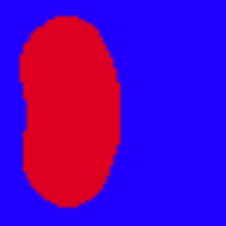} } \vskip2pt
\centerline{
\includegraphics[width=0.24\textwidth]{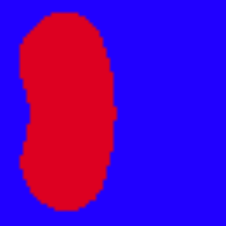}
\includegraphics[width=0.24\textwidth]{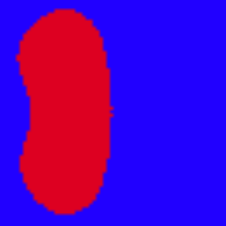}
\includegraphics[width=0.24\textwidth]{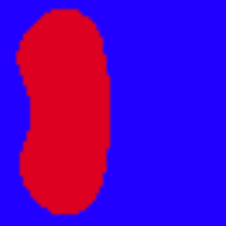}
\includegraphics[width=0.24\textwidth]{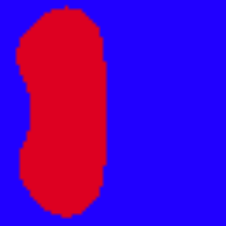} } \vskip2pt
\centerline{
\includegraphics[width=0.24\textwidth]{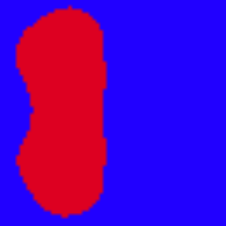}
\includegraphics[width=0.24\textwidth]{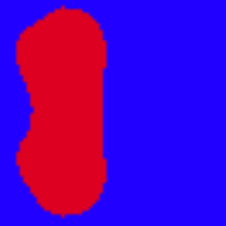}
\includegraphics[width=0.24\textwidth]{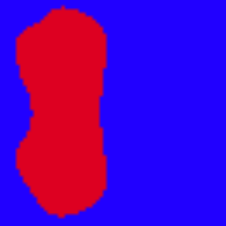}
\includegraphics[width=0.24\textwidth]{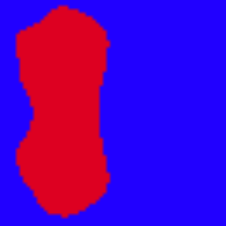} }
\caption{Level set evolution for second experiment. Plots after
$0,1,5,20$, $50,100,200,400$, $600,800,1000,2000$, $5000,10000,20000,50000$
iterative steps. \label{fig:ls-evol-4} }
\end{figure}

As in the first part of the previous experiment, the data is
almost given exactly and the velocity $w_k$ is again obtained by
solving the boundary value problem with $\beta=0$.
The evolution of the level set function is
shown in Figure~\ref{fig:ls-evol-4}.

\begin{remark}
The effect of parameter changes: In our numerical observations  we observed that in numerical simulations the minimizer is not
severely affected by the choice of $\beta \alpha$ and can in fact
be neglected.
\end{remark}


\section*{Acknowledgment}
The work of F.F. has been supported by the Tiroler
Zukunftsstiftung; the work of O.S. has been partly supported by
the FWF (Austrian Science Foundation), grant Y-123 INF.
A.L. is on leave from Department of Mathematics, Federal Univ. of
St.\,Catarina, Brazil; his work is supported by the Austrian Academy
of Sciences and CNPq, grant 305823/2003-5.

\end{document}